\documentclass[12pt,a4paper,titlepage,fleqn]{article}
 
\usepackage{amssymb,amsmath,stmaryrd}
 
\newtheorem{theorem}{Theorem}[section]
\newtheorem{lemma}{Lemma}[section]

\newtheorem{corollary}{Corollary}[section]
\newtheorem{proposition}{Proposition}[section]

\newcommand{\sgn}{{\rm sgn}}
\newcommand{\divergence}{{\rm div}}

\def\qed{\hfill$\square$\\ $\phantom{A}$\\} 
\def\barint{\kern4pt \raise3.4pt\hbox{\vrule height.8pt width6pt} \kern-10pt 
\int}

\begin{document}
\label{begindoc} 

\title
{Coercive Inequalities on Metric Measure Spaces.
\thanks{
{Supported by EPSRC 
EP/D05379X/1} \& Royal Society.  
The first author was partially supported by 
KBN grant {\tt 1 P03A 03029}} 
}
\author{W. Hebisch \\
{\small{Institute of Mathematics}}\\
{\small{University of Wroc{\l}aw}}\\
{\small{Wroc{\l}aw, Poland
}}\\
\\ $\ $\\
B. Zegarli{\'n}ski\\
{\small{Department of Mathematics}}\\
{\small{Imperial College}}\\
{\small{London, UK}}
}

\date{ }
\maketitle

 
 
 
 
\newpage

\section{Introduction} 
\label{1.Intro}

In this paper we study coercive inequalities on finite dimensional metric spaces with probability measures which do not have volume doubling property. This class of inequalities includes the well known Poincar\'e inequality
$$
M \mu |f -\mu f|^q \leq \mu |\nabla f|^q
$$
with some constants $M\in(0,\infty)$, $q\in(1,\infty)$ independent a function $f$ for which the (metric) length of the gradient $|\nabla f|$ is well defined,
as well as a variety of stronger coercive inequalities with the variance on the left hand side replaced by a functional with a stronger growth, as for example in case of celebrated Log - Sobolev inequality which is of the following form 
$$
\mu f^2 \log \frac{f^2}{\mu f^2} \leq c\mu |\nabla f|^2
$$
with some constant $c\in(0,\infty)$ independent of a function $f$.

We are interested in probability measures on noncompact spaces, like for example the finite products of real lines $\mathbb{R}^n$, but also certain noncompact groups as for example the Heisenberg group. 

For probability measures on the real line the necessary and sufficient condition for Poincar\'e inequality characterising the density (of the absolutely continuous part with respect to the Lebesgue measure) were established long time ago by Muckenphout, \cite{Mu}, (\cite{Ma}). More recently such criteria were established for 
other coercive inequalities (Log-Sobolev type: (LS$_2$) \cite{BG} , (LS$_q$) \cite{BZ}, for distributions with weaker tails \cite{BCR},...). In multidimensional case the situation is rather different and more intricate. First of all, since the inequalities of interest to us have a natural tensorisation property, there is a number of perturbative techniques which allow to obtain classes of interesting examples in higher and even in infinite dimensions (see e.g. \cite{GZ}, \cite{BZ}, \cite{RZ}, \cite{LZ},.., \cite{BoHe},\cite{Y},.. and references given there). 
We would like to mention a work \cite{R} in which the coercive inequalities for probability measures on $\mathbb{R}^n$, $n\geq 3$, with variety of decay of the tails
(slower as well as faster than the Gaussian) were systematically studied with the help of classical Sobolev inequalities providing in particular an effective
sufficient criteria, (in terms of certain nonlinear differential inequalities 
for the log of the density function), for related coercive inequalities,
(see also reviews \cite{Si}, \cite{DaGrSi} and references therein).
In the mid 80'ties Bakry and Emery, \cite{Ba}, introduced a very effective criterion based on convexity (curvature) which allowed to enlarge a class of examples where Log-Sobolev inequality holds, including situation with measures on certain finite dimensional Riemannian manifolds; (as well as some infinite dimensional cases however with a compact configuration space \cite{CS}).
Following similar line of reasoning, in \cite{AA}  the authors provided an
effective criteria for (generalisation) of Brascamp-Lieb inequality as well as Log-Sobolev inequality (with possibly more general entropy functional and weighted
Dirichlet form dependent on the measure).

More recently, in \cite{BL1}, certain convexity ideas, (including  Brunn-Minkowski inequality), 
were exploited to recover in the special case of the space $\mathbb{R}^n$ similar results as in \cite{AA} and obtained additionally inequalities (LS$_q$) which are naturally related to different than Euclidean metrics (in particular involving different length of the gradient on the right hand side). These result concerned principally the probability measures with tails decaying faster than the Gaussian.
We point out that while such distribution were also discussed in \cite{R}, in \cite{BL1} they involved in a natural way Lipschitz functions with respect to a non-Euclidean metric (while in Rosen's work the emphasis of improvement was on different functionals on the left hand side). The corresponding results for measures on $\mathbb{R}^n$ with slower distribution tails were obtained in \cite{BCR} (see also references therein), which
included in particular those of Rosen, \cite{R}, for the similar class of measures.

Part of the motivation for the current paper was provided by \cite{LZ} in which
the coercive inequalities involving H\"ormander fields instead of the (nondegenerate full gradient) were studied. Such the situation is naturally related to a more general Carnot-Caratheodory metric associated to the family of fields and the interest here is to obtain coercive inequalities involving length of the corresponding metric gradient. While in \cite{LZ} a rich family of examples on compact spaces was provided, the noncompact situation was more difficult. In this paper we develop an efficient technology which not only recovers interesting results in $\mathbb{R}^n$ briefly reviewed in the above, but also allows us to extend to interesting metric spaces
as certain noncompact Lie groups including in particular the Heisenberg group.
Part of our approach is directed on proving inequalities, which we call U-bounds, 
of the following form
$$
\int |f|^q U d\mu \leq C \int |\nabla f|^q d\mu + D \int |f|^q d\mu
$$
with a suitable increasing unbounded function $U$ of the metric and the length of the metric gradient $|\nabla f|$; see Section 2. We show later in Section 3 and 4 that
such inequality implies corresponding Poincar\'e as well as suitable coercive inequalities; in fact as we illustrate in some of the cases the U-bounds are equivalent with the coercive inequalities.
(This requires an extension of result of on a Gaussian exponential bound of \cite{AiSt} for other measures and functions with possibly unbounded gradient.)\\
In Section 5 we explore also a family of weighted Poincar\'e and Log-Sobolev inequalities on Riemannian manifolds including measures with ultra slow tails. In such the context we can effectively employ Laplacian comparison theorem (see e.g. \cite{CLN}) which in particular allows us to extend recent results of \cite{BL2} where convexity ideas in Euclidean spaces were used.
\\
As an application of our technique we also prove (see Section 6-7) the Log-Sobolev inequality
for the heat kernel measure on the Heisenberg group, (a topic which attracted recently some extra attention \cite{Li}, \cite{DM}). 



\section{ U - Bounds.}
\label{2.Operator Bounds}

By $\nabla$ we denote a subgradient in $\mathbb{R}^N$, that is a finite collection of possibly noncommuting fields. It is assumed that the divergence of each of these fields with respect to the Lebesgue measure
$\Lambda$ on $\mathbb{R}^N$ is zero. (While this provides some simplification in our expositions, it is possible to extend our arguments to a more general setting.)\\
We begin with proving the following result.

\begin{theorem}\label{thm2.5}
Let $d\mu_p = \frac{e^{-\beta d^p}}{Z} d\lambda$ be a probability
measure defined with $\beta\in(0,\infty)$ and $p\in(1,\infty)$,
($Z$ being the normalisation constant). 
Suppose $0<\frac1\sigma \leq |\nabla d| 
\leq 1 $, for some $\sigma\in[1,\infty)$, and
$\Delta d\leq  K+ \beta p \varepsilon d^{p-1}$ outside the unit ball $B\equiv \{d(x)< 1\}$ 
for some $K\in[0,\infty)$ and $\varepsilon \in[0,\frac1{\sigma^2})$.
Then there exist constants $C,D\in(0,\infty)$ such that
the following bound is true.
\begin{equation} \label{2.z}
\int |f|d^{p-1}\,  d\mu_p \leq C \int |\nabla f|\, d\mu_p + D \int |f|\, d\mu_p 
\end{equation}
\end{theorem}

\noindent\textbf{Remark}: In particular the assumptions of the theorem are satisfied for $d$ being the Carnot-Caratheodory
distance and $\nabla$ the (horizontal) gradient of the Heisenberg group.\\

 \noindent\textit{Proof}:
For a smooth function $f \geq 0$ such that $f = 0$ on the
unit ball, by the Leibniz rule we have
\begin{equation} \label{2.10}
(\nabla f) e^{-\beta d^p}= \nabla \left(f  e^{-\beta d^p}\right)
 + \beta pf\left( d^{p-1}\nabla d\right)e^{-\beta d^p}.
\end{equation}
Put
$$\boldsymbol{\alpha}(\cdot)\equiv \int (\nabla d) ( \cdot )\  d\lambda.$$
Acting with this functional on the expression (\ref{2.10})  we get
\begin{equation} \label{2.11}
\boldsymbol{\alpha}((\nabla f) e^{-\beta d^p}) = \boldsymbol{\alpha}
\left(\nabla \left(f  e^{-\beta d^p}\right)\right) +
\beta p\int f d^{p-1} |\nabla d|^2\  e^{- \beta d^p}d\lambda
\end{equation}
Using H\"older inequality, the left hand side of (\ref{2.11}) can be
estimated from above as follows
\begin{eqnarray} \label{2.12}
\boldsymbol{\alpha}((\nabla f) e^{-\beta d^p}) =
\int (\nabla d) \cdot (\nabla f) e^{-\beta d^p}d\lambda \phantom{AAAAAAAAA} \\
\phantom{AAAAAAAAa}\leq \int |\nabla d||\nabla f|e^{-\beta d^p}d\lambda 
\leq \int |\nabla f|e^{-\beta d^p}d\lambda \nonumber
\end{eqnarray}
where we have used the fact that $|\nabla d |\leq 1$.
The first term on the right hand side of (\ref{2.11}) can be
treated with the help of integration by parts as follows
\begin{eqnarray} \label{2.13}
\boldsymbol{\alpha}\left(\nabla \left(f  e^{-\beta d^p}\right)\right)
=
\int (\nabla d) \cdot \nabla \left(f  e^{-\beta d^p}\right)d\lambda
 \\
= - \int (\Delta d) f e^{-\beta d^p}d\lambda 
\geq -K\int f e^{-\beta d^p}d\lambda  -   \beta p\varepsilon \int f d^{p-1} e^{-\beta d^p}d\lambda
\nonumber 
\end{eqnarray}
where we have used the assumption that $\Delta d \leq K+\beta p \varepsilon d^{p-1}$.
Combining (\ref{2.11}), (\ref{2.12}) and (\ref{2.13}), we get
$$
\beta p\int f d^{p-1} \left(|\nabla d|^2-\varepsilon\right)\  e^{- \beta d^p}d\lambda \leq
\int |\nabla f|e^{-\beta d^p}d\lambda + K\int f e^{-\beta d^p}d\lambda
$$
from which the inequality  (\ref{2.z}) follows
with $C=\frac{1}{(1/\sigma^2-\varepsilon)\beta p}$ and $D=\frac{K}{(1/\sigma^2-\varepsilon)\beta p}$,
provided $\varepsilon \in[0,\frac1{\sigma^2})$.

Now, the estimate (\ref{2.z}) is proven for smooth nonnegative $f$ which
vanish on the unit ball.  We can handle non-smooth functions 
approximating them by
smooth ones (on compact sets via convolution and splitting f into
compactly supported pieces using a smooth partition of unity --
details are tedious but do not pose any essential difficulty).

We can handle $f$ of arbitrary sign replacing $f$ by
$|f|$ and using equality $\nabla |f| = \sgn(f)\nabla f$.

To handle $f$ which are non-zero on the unit ball
we write $f = f_0 + f_1$ where $f_0 = \phi f$, $f_1 = (1-\phi)f$ and
$\phi(x) = \min(1, \max(2 - d(x),0))$.  Then
$$
\int |f| d^{p-1}d\mu_p = \int_{d(x)\leq 2}|f| d^{p-1}d\mu_p 
+ \int_{d(x)>2}|f| d^{p-1}d\mu_p$$
$$
\leq 2^{p-1}\int_{d(x)\leq 2}|f|d\mu_p + \int_{d(x)>2}|f|_1 d^{p-1}d\mu_p$$
$$
\leq 2^{p-1}\int |f|d\mu_p + \int |f|_1 d^{p-1}d\mu_p.$$
Next
$$|\nabla f_1| \leq |\nabla f|+|f|,$$
$$\int |f_1| d^{p-1}d\mu_p
\leq C\int |\nabla f_1|d\mu_p + D\int |f_1|d\mu_p$$
$$
\leq C\int |\nabla f|d\mu_p + (D+C)\int |f|d\mu_p$$
Combining inequalities above we see that (\ref{2.z}) is valid without
restriction on the support of $f$ if we replace $D$ by $D+2^{p-1} + C$.
\qed

Using our result and a perturbation technique we obtain the following
generalisation.
\begin{theorem}\label{thm2.6}
Let $d\mu = \frac{e^{-W-V}}{Z'} d\mu_\theta$ be a probability measure defined with a differentiable potential $W$ satisfying
\begin{equation}\label{2x.7c}
|\nabla W | \leq \delta d^{p-1} + \gamma_\delta
\end{equation}
with some constants $\delta < 1/C$ and $\gamma_\delta\in(0,\infty)$,
and suppose that $V$ is a measurable function such that $osc(V)\equiv \max V - \min V <\infty $.
Then there exist constants $C',D'\in(0,\infty)$ such that
the following bound is true.
\begin{equation}\label{2x.8c}
\int |f| d^{p-1} d\mu \leq C' \int |\nabla f| d\mu + D' \int |f| d\mu
\end{equation}
\end{theorem}

\noindent\textbf{Remark}: In particular the assumption (\ref{2x.7c}) of the
theorem is satisfied if $W$ is a polynomial of lower order
in $d$. Another example, in the spirit of \cite{GR} and \cite{BZ},
with deep wells is as follows
$$
W = \vartheta d^{p-1} \cos(d)
$$
with a small constant $\vartheta>0$, (but 
$\vartheta d^{p-1} \cos(d^{1+\varepsilon})$ would not work for any
$\varepsilon>0$ no matter how small $\vartheta>0$ would be).

\medskip

\textit{Proof}:
We consider first the case $V=0$ and start from substituting
$fe^{-W}$ in the inequality (\ref{2.z}) for the measure
$\mu_p$. Using Leibniz rule
$$
\int |f| d^{p-1} e^{- W} d\mu_p \leq C\int |\nabla f| e^{-W} d\mu_p + D
\int |f| e^{-W} d\mu_p
$$
$$\phantom{AAAAAAAAAAAAAAAAA}+
C\int |f||\nabla W| e^{- W} d\mu_p
$$
Now our assumption (\ref{2x.7c}) about $W$, implies
$$
\int |f||\nabla W| e^{-W}\, d\mu_p \leq 
\delta \int |f| d^{p-1}  e^{-W}\, d\mu_p  +
\gamma_\delta \int |f| e^{-W}\, d\mu_p
$$
Thus combining these bounds we arrive at
$$
\int |f| d^{p-1} e^{-W} d\, \mu_p \leq 
\bar C\int |\nabla f| e^{- W}\, d\mu_p +
\bar D \int |f| e^{- W}\, d\mu_p
$$
with
$$ \bar C\equiv C / (1- C\delta)
\quad \textrm{and} \quad
\bar D = (D + \gamma_\delta)/ (1- C\delta)
$$
Next we note that if $V\neq 0$ we have
\begin{eqnarray*}
\int |f| d^{p-1} \frac {e^{- W-V}}{Z'}\, d\mu_p
\leq  e^{osc(V)}\int |f| d^{p-1} \frac {e^{- W}}{\int e^{- W}\, d\mu_p}\, d\mu_p
\phantom{AAAAAAAAAAAAAAAAAAAA}\\
\\
\leq e^{osc(V)} \bar C \int |\nabla f|  \frac {e^{- W}}{\int e^{- W} d\mu_p}\,
 d\mu_p
+ e^{osc(V)} \bar D \int |f|  \frac {e^{- W}}{\int e^{- W} d\mu_p}\, d\mu_p
\phantom{AAAAAAAA}\\
\\
\leq   e^{2osc(V)} \bar C \int |\nabla f| \frac {e^{- W-V}}{Z'}\,  d\mu_p
+ e^{2osc(V)} \bar D \int |f|  \frac {e^{- W-V}}{Z'}\, d\mu_p
\phantom{AAAAAAAAAA}
\end{eqnarray*}
\qed

\begin{theorem}\label{thm2.5p}
Let $\mu$ be a probability measure for which conclusion of
Theorem \ref{thm2.5} holds.
Let $p\in(1,\infty)$.
Then for each $q\in[1,\infty)$
there exist constants $C_q,D_q\in(0,\infty)$ such that
the following bound is true.
\begin{equation} \label{2.ap}
\int |f|^q d^{q(p-1)}\,  d\mu 
\leq C_q \int |\nabla f|^q\, d\mu + D_q \int |f|^q\, d\mu.
\end{equation}
\end{theorem}

\textit{Proof}: Let $d_1(x) = \max(1, d(x))$.  Enlarging constants
$D$ if necessary we may assume that
$$
\int |f|d_1^{p-1}\, d\mu
\leq C \int |\nabla f|\,  d\mu + D\int |f|\,  d\mu.$$
Put $h = |f|^qd_1^{(p-1)(q-1)}$.  We have
$$
\int |f|^q d^{q(p-1)}\,  d\mu 
\leq \int |f|^q d_1^{q(p-1)}\,  d\mu = \int h d_1^{p-1}\,  d\mu$$
$$
\leq
C \int |\nabla h|\,  d\mu + D\int h\,  d\mu.$$
By Leibniz formula
$$
|\nabla h| = q|\nabla f||f|^{(q-1)}d_1^{(q-1)(p-1)} 
+ (q-1)(p-1)|\nabla d_1| |f|^q d_1^{(q-1)(p-1)-1}$$
and
$$
\int q|\nabla f||f|^{(q-1)}d_1^{(q-1)(p-1)} \,  d\mu$$
$$
\leq q\left(\int |\nabla f|^qd\mu\right)^{1/q}
\left(\int(|f|^{q-1}d_1^{(q-1)(p-1)})^{q/(q-1)} d\mu \right)^{(q-1)/q}
$$
$$
\leq \alpha^q\int |\nabla f|^qd\mu + \frac{q-1}{\alpha^{q/(q-1)}}
\int |f|^{q}d_1^{q(p-1)}d\mu.
$$
Next
$$
\int h\,  d\mu = \int |f|^qd_1^{(q-1)(p-1)}\,  d\mu \leq
\left(\int |f|^q\,  d\mu\right)^{1/q}
\left(\int|f|^q d_1^{q(p-1)}\,  d\mu \right)^{(q-1)/q}
$$
$$
\leq \frac{\beta^q}{q}\int |f|^q\,  d\mu + 
\frac{q-1}{\beta^{q/(q-1)}q}\int|f|^q d_1^{q(p-1)}\,  d\mu.
$$
If $(q-1)(p-1) \leq 1$, then
$$
\int (q-1)(p-1)|\nabla d_1| |f|^q d_1^{(q-1)(p-1)-1}\, d\mu \leq
\int |f|^q\, d\mu.$$
If $(q-1)(p-1) > 1$, then
$$
\int (q-1)(p-1)|\nabla d_1| |f|^q d_1^{(q-1)(p-1)-1}\, d\mu
\leq (q-1)(p-1) \int |f|^q d_1^{(q-1)(p-1)-1}\, d\mu $$
$$
\leq (q-1)(p-1) \left(\int |f|^q\, d\mu\right)^{p/(q(p-1))}
\left(\int  |f|^qd_1^{q(p-1)}\, d\mu\right)^{((p-1)(q-1)-1)/(q(p-1))}$$
$$
\leq \frac{(q-1)p}{q}\gamma^{q(p-1)/p}\int |f|^q\, d\mu +
\frac{(q-1)^2(p-1)-(q-1)}{q\gamma^{q(p-1)/((q-1)(p-1)-1)}}\int  |f|^qd_1^{q(p-1)}\, d\mu
$$
Combining inequalities above, if $(q-1)(p-1) \leq 1$ we get 
$$ 
(1 - C \frac{q-1}{\alpha^{q/(q-1)}} -D\frac{q-1}{\beta^{q/(q-1)}q})
\int  |f|^qd_1^{q(p-1)}\, d\mu
$$
$$
\leq C \alpha^q\int |\nabla f|^qd\mu 
+ (C +  D\frac{\beta^q}{q})\int |f|^q\, d\mu $$
which gives the claim with
$$
C_q = \frac{C\alpha^q}
{1 - C \frac{q-1}{\alpha^{q/(q-1)}} -D\frac{q-1}{\beta^{q/(q-1)}q}}$$
and
$$
D_q = \frac{C +  D\frac{\beta^q}{q}}
{1 - C \frac{q-1}{\alpha^{q/(q-1)}} -D\frac{q-1}{\beta^{q/(q-1)}q}}$$
if $\alpha$ and $\beta$ are big enough.  Similarly, for
$(q-1)(p-1) > 1$ we get the claim with
$$
C_q = \frac {C\alpha^q}
{1 - C \frac{q-1}{\alpha^{q/(q-1)}} -C\frac{(q-1)^2(p-1)-(q-1)}
{q\gamma^{q(p-1)/((q-1)(p-1)-1)}}
-D\frac{q-1}{\beta^{q/(q-1)}q}}
$$
and 
$$
D_q = \frac{Cp\gamma^{q(p-1)/p} +  D\frac{\beta^q}{q}}
{1 - C \frac{q-1}{\alpha^{q/(q-1)}} -C\frac{(q-1)^2(p-1)-(q-1)}
{q\gamma^{q(p-1)/((q-1)(p-1)-1)}}
-D\frac{q-1}{\beta^{q/(q-1)}q}}
$$
if $\alpha$, $\beta$ and $\gamma$ are big enough.
\qed


\begin{theorem}\label{thm2.1}
Let $d\mu_p = \frac{e^{-\beta d^p}}{Z} d\lambda$ be a probability measure defined with $\beta\in(0,\infty)$ and $p\in[2,\infty)$,
($Z$ being the normalisation constant). Suppose $0<\frac1\sigma \leq |\nabla d| \leq 1 $, for some $\sigma\in[1,\infty)$, and
$\Delta d\leq  K+ \beta p \varepsilon d^{p-1}$ outside the unit ball $B\equiv \{d(x)< 1\}$ 
for some $K\in[0,\infty)$ and $\varepsilon \in[0,\frac1{\sigma^2})$.

\noindent Suppose $\frac1q+\frac1p=1$, then we have
\begin{equation}\label{2.ap_bis}
\int |f|^q d^p d\mu \leq C_q \int |\nabla f|^q d\mu + D_q \int |f|^q d\mu
\end{equation}
\end{theorem}

\noindent\textbf{Remark}: In particular the assumptions of the theorem are satisfied for $d$ being the Carnot-Caratheodory
distance and $\nabla$ the (horizontal) gradient of the Heisenberg group.

\textit{Proof}: This is a special case of Theorem \ref {thm2.5p}
\qed

Extension to more general measures is as follows.
\begin{theorem}\label{thm2.2}
Let $d\mu = \frac{e^{-W-V}}{Z'} d\mu_p$ be a probability measure defined with a differentiable potential $W$ satisfying 
\begin{equation}\label{2.7a}
|\nabla W |^q \leq \delta d^p + \gamma_\delta
\end{equation}
with some constants $\delta 2^{q-1}q^{-q}C < 1$ and $\gamma_\delta\in(0,\infty)$,
and suppose that $V$ is a measurable function such that $osc(V)\equiv \max V - \min V <\infty $.
Then there exist constants $C',D'\in(0,\infty)$ such that
the following bound is true.
\begin{equation}\label{2.8a}
\int |f|^q d^p d\mu \leq C' \int |\nabla f|^q d\mu + D' \int |f|^q d\mu 
\end{equation}
with $q$ such that $\frac1q+\frac1p =1$.
\end{theorem}

The proof is similar to that of Theorem \ref{thm2.6}

\subsection{U - Bounds: Sub-Quadratic Case.}%
\label{2x.U Bounds Bis}

\begin{theorem}\label{thm2.3}
Let $d\mu_\theta = \frac{e^{-\beta d^\theta}}{Z} d\lambda$ be a probability measure defined with $\beta\in(0,\infty)$ and $\theta\in[1,2)$,
($Z$ being a normalisation constant). Suppose $0<\frac1\sigma \leq |\nabla d| \leq 1 $, for some $\sigma\in[1,\infty)$, and
$\Delta d\leq  K+ \beta p \varepsilon d^{p-1}$ outside the unit ball $B\equiv \{d(x)< 1\}$ 
for some $K\in[0,\infty)$ and $\varepsilon \in[0,\frac1{\sigma^2})$.\\
Then there exist constants $C_\theta,D_\theta\in(0,\infty)$ such that
the following bound is true
\begin{equation} \label{2x.a}
\int |f|^2 d^{2(\theta-1)} d\mu_\theta \leq 
C_\theta \int |\nabla f|^2 d\mu_\theta + D_\theta \int |f|^2 d\mu_\theta 
\end{equation}
\end{theorem}

\noindent\textbf{Remark}: In particular the assumptions of the theorem are satisfied for $d$ being the Carnot-Caratheodory
distance and $\nabla$ the (horizontal) gradient of the Heisenberg group.

\noindent\noindent\textit{Proof}: Again, this is a special
case of Theorem \ref{thm2.5p}
\qed

Extension to more general measures is as follows.
\begin{theorem}\label{thm2.4}
Let $d\mu = \frac{e^{-W-V}}{Z'} d\mu_\theta$ be a probability measure defined with a differentiable potential $W$ satisfying 
\begin{equation}\label{2x.7a}
|\nabla W |^2 \leq \delta d^{2(\theta-1)} + \gamma_\delta
\end{equation}
with some constants $\delta C/2 < 1$ and $\gamma_\delta\in(0,\infty)$,
and suppose that $V$ is a measurable function such that $osc(V)\equiv \max V - \min V <\infty $.
Then there exist constants $C',D'\in(0,\infty)$ such that
the following bound is true.
\begin{equation}\label{2x.8a}
\int |f|^2 d^\theta d\mu \leq C' \int |\nabla f|^2 d\mu + D' \int |f|^2 d\mu 
\end{equation}
\end{theorem}
Again, the proof is similar to that of Theorem \ref{thm2.6}

\section{Poincar\'e inequality.}
\label{3.PoicareIneq}

\begin{theorem} \label{thm3.1}
Suppose $1 \leq q < \infty$ and a measure $\lambda$ satisfies the $q$-Poincar\'e inequality for every ball $B_R$, that is there 
exists a constant $c_R\in(0,\infty)$ such that
\begin{equation} \label{eq3.0a}
\frac1{|B_R|}\int_{B_R} \left|f- \frac1{|B_R|}\int_{B_R}f\right|^q d\lambda \leq c_R \frac1{|B_R|}\int_{B_R} |\nabla f|^q d\lambda
\end{equation}
Let $\mu$ be a probability measure on $\mathbb{R}^n$ which is absolutely continuous
with respect to the measure $\lambda$ and such that
\begin{equation} \label{eq3.0b}
\int f^q \eta d\mu \leq C \int |\nabla f|^q d\mu + D \int f^q d\mu
\end{equation}
with some nonnegative function $\eta$ and some constants $C,D\in(0,\infty)$ independent of a function $f$. 
If for any $L\in(0,\infty)$ there is a constant $A_L$ such that 
\begin{equation} \label{eq3.0c}
\frac1{A_L} \leq \frac {d\mu}{d\lambda} \leq A_L
\end{equation}
on the set $\{\eta < L\}$ and, for some $R\in(0,\infty)$ (depending on $L$),
 we have $\{\eta < L\}\subset B_R$,
then $\mu$ satisfies the $q$-Poincar\'e inequality
\begin{equation}\label{eq3.0d}
\mu \left|f- \mu f\right|^q \leq c \mu |\nabla f|^q
\end{equation}
\end{theorem}

\textit{Proof}: For any $a$ we have 
\begin{equation} \label{eq3.1prim}
\mu \left|f- \mu f\right|^q \leq 2^q \mu \left|f- a\right|^q.
\phantom{AAAAAAAAAAAAAAAAAAAAAAAAAAAAA}
\end{equation}
Next
\begin{equation} \label{eq3.1}
\mu \left|f- a\right|^q
 \leq \mu  \left|f- a\right|^q \chi(\eta< L) +
 \mu \left|f- a\right|^q \chi(\eta\geq L) 
\end{equation}
Using our assumptions and putting $a = \frac1{|B_R|}\int_{B_R}f$, 
for the first term on the right hand side of (\ref{eq3.1}) we have
\begin{eqnarray} \label{eq3.2}
 \mu  \left|f- a\right|^q \chi(\eta < L)
\leq 
 A_L\int_{B_R} \left|f- \frac1{|B_R|}\int_{B_R}f\right|^q d\lambda
\nonumber \\
\leq A_L c_R  
\int_{B_R} |\nabla f|^q d\lambda  \leq  A_L^2 c_R \mu |\nabla f|^q 
\end{eqnarray}
On the other hand for the second term on the right hand side of (\ref{eq3.1}) we get
\begin{eqnarray} \label{eq3.3}
\mu \left|f- a\right|^q \chi(\eta \geq L) 
\leq \frac{1}{L}\mu \left|f-a \right|^q \eta 
\end{eqnarray}
Hence, by (\ref{eq3.0b}), we obtain
\begin{equation} \label{eq3.4}
\mu \left|f- a\right|^q \chi(\eta \geq L) 
\leq  \frac{C}{L} \mu|\nabla f|^q+ \frac{D}{L} \mu\left|f- a\right|^q
\end{equation}
Combining (\ref{eq3.2}) and (\ref{eq3.4}), we get
$$
\mu \left|f- a\right|^q \leq \left[A_L^2  c_R +  \frac{C}{L} \right] \mu |\nabla f|^2  +
\frac{D}{L} \mu\left|f- a \right|^q
$$
Choosing $L > D $,
simple rearrangement yields
$$
\mu \left|f - a\right|^q \leq \frac{A_L^2 c_R +  \frac{C}{L}}{1-\frac{D}{L}}\  \mu |\nabla f|^q 
$$
This together with (\ref{eq3.1prim}) - (\ref{eq3.2}) yields
$$
\mu |f-\mu f|^q \leq c \mu |\nabla f|^q
$$
with some constant $c\in(0,\infty)$.
\qed
\begin{corollary} \label{thm3.2}
If we are on nilpotent Lie group the probability measure $\mu_q$ and $\mu_\theta$ of Theorem \ref{thm2.2} and \ref{thm2.4}, respectively,
satisfies the Poincar\'e inequality. 
\end{corollary}

\section{From Sobolev Inequalities to Coercive Inequalities with Probability Measure: 
The non-compact setting.} 
\label{4.Non-compact setting}

\subsection{Case $p\geq 2$.} 

\begin{theorem} \label{thm4.1}
Let $d\mu= \frac{e^{-U}}{Z}d\lambda$.
Suppose the following Sobolev inequality is satisfied 
\begin{equation} \label{eq4.0a}
\quad \left(\int |f|^{q+\varepsilon} d\lambda\right)^{\frac{q}{q+\varepsilon}} \leq
a \int |\nabla f|^q  d\lambda + b\int |f|^q  d\lambda
\end{equation}
and the following bound is true
\begin{equation}\label{eq4.0b}
\mu \left(|f|^q \left[|\nabla U|^q + U\right] \right) \leq \bar C\mu|\nabla f|^q + \bar D \mu |f|^q
\end{equation}
Then the following inequality is true
\begin{equation}\label{eq4.0c}
\quad\,\, \mu \left(f^q \log \frac {f^q}{\mu f^q}\right) \leq C \mu|\nabla f|^q + D \mu |f|^q 
\end{equation}
Moreover, if $q\in(1,2]$ and the following q-Poincar\'e inequality holds
\begin{equation}\label{eq4.0d}
\qquad \qquad   \mu|f-\mu f|^q \leq \frac1M \mu|\nabla f|^q
\end{equation}
then one has
\begin{equation}\label{eq4.0e}
\quad \, \mu \left(f^q \log \frac {f^q}{\mu f^q}\right) \leq c \mu|\nabla f|^q
\end{equation}
with some constant $c\in(0,\infty)$ independent of $f$. 
\end{theorem}

\textit{Proof}: First we note that for $f\not\equiv 0$, we have
$$
\mu \left(f^q \log \frac {f^q}{\mu f^q}\right) = \mu (f^q) \int g^q \log {g^q} d\lambda + 
\mu \left(f^q [U+\log Z]\right) 
$$
with $g\equiv f\cdot \frac{e^{-\frac1q U}}{Z^\frac1q}$ satisfying  $\int g^q d\lambda =1$.
Next, by arguments based on Jensen inequality, one gets
$$
\int g^q \log {g^q} d\lambda = \frac q\varepsilon  \int g^q \log {g^\varepsilon} d\lambda \leq
\frac {q(q+\varepsilon)}{\varepsilon}  \log \left(\int g^{q+\varepsilon} d\lambda\right)^{\frac1{q+\varepsilon}}
$$
whence, by the Sobolev inequality (\ref{eq4.0a}), one obtains
$$
\int g^q \log {g^q} d\lambda\leq \frac {q+\varepsilon}{\varepsilon}   \log \left(\int g^{q+\varepsilon} d\lambda\right)^{\frac q{q+\varepsilon}} 
\leq a' \int |\nabla g|^q d\lambda + b' \int g^q d\lambda
$$
with $a'\equiv \frac {q+\varepsilon}{\varepsilon} a$ and $b'\equiv\frac {q+\varepsilon}{\varepsilon}b $.
Combining all the above we arrive at
$$
\mu \left(f^q \log \frac {f^q}{\mu f^q}\right) \leq
a'\mu(f^q) \int |\nabla \left(f\frac{e^{-\frac1q U}}{Z^\frac1q}\right) |^q d\lambda + 
(b'+\log Z) \int f^q d\mu
+ \mu \left(f^q U\right)
$$
and, by simple arguments, we obtain
\begin{equation} \label{eq4.1}
 \mu \left(f^q \log \frac {f^q}{\mu f^q}\right) \leq
2^{q-1}a' \int |\nabla f|^q d\mu 
+ \mu \left(f^q \left[2^{q-1}q^{-q}a' |\nabla U|^q + U + b' +\log Z \right] \right)
\end{equation}
Now using our assumption (\ref{eq4.0b}) yields
\begin{equation} \label{eq4.2}
\mu \left(f^q \log \frac {f^q}{\mu f^q}\right) \leq \left( 2^{q-1}a'  +  2^{q-1}q^{-q}a' \bar C \right)
\mu|\nabla f|^q + (b'+\bar D+\log Z) \mu |f|^q
\end{equation}
Since for $q\in(1,2]$ one has, \cite{BZ},
\begin{equation}\label{4.*}
\mu \left(f^q \log \frac {f^q}{\mu f^q}\right) \leq 
\mu \left(|f-\mu f|^q \log \frac {|f-\mu f|^q}{\mu |f-\mu f|^q}\right) + 2^{q+1}\mu |f-\mu f|^q
\end{equation}
using (\ref{eq4.2}) we arrive at
$$
\mu \left(f^q \log \frac {f^q}{\mu f^q}\right) \leq 
\left\{\left( 2^{q-1}a'  +  2^{q-1}q^{-q}a' \bar C \right) +\frac{2^{q+1}(b'+\bar D+\log Z)}{M}
\right\} \mu |\nabla f|^q
$$
which ends the proof of the theorem.
\qed 
Using Theorem \ref{thm4.1} together with results of Section 3, (q-Poincar\'e inequality), we arrive at
the following result.

\begin{corollary} \label{cor4.2}
The probability measures $d\mu = e^{-W-V}d\mu_p/ Z'$, with $p\geq 2$, described in Theorem \ref{thm2.2}
satisfies the following coercive inequality
\begin{equation*}
\phantom{AAAAAAAAAAAA}
\mu\left(|f|^q\log \frac{|f|^q}{\mu |f|^q}\right) \leq c \mu |\nabla f|^q
\tag{\textbf{LS}$_q$}
\end{equation*}
with $\frac1q+\frac1p=1$ and a constant $c\in(0,\infty)$ independent of a function $f$.
\end{corollary}

\subsection{Sub-quadratic Case.} 

\begin{theorem} \label{thm4.2}
Suppose $\theta\in[1,2]$ and let 
$\varsigma = \frac{2(\theta-1)}{\theta}$. Then there exist constants $C,D\in(0,\infty)$ such that
\begin{equation} \label{4.2.a}
\int f^2 \left|\log \frac{f^2}{\int f^2d\mu_\theta}\right|^\varsigma d\mu_\theta 
\leq C  \int |\nabla f|^2 d\mu_\theta + D \int f^2 d\mu_\theta
\end{equation}

\end{theorem}

\textit{Proof}: 
We note first that if $\theta\in[1,2]$, then $\varsigma \in [0,1]$.
Put $g\equiv f \frac{e^{-\frac\beta2 d^\theta}}{Z^{\frac12}}$.
We have the following inequality
\begin{eqnarray} \label{4.2.1}
\int f^2 \left|\log \frac{f^2}{\int f^2d\mu_\theta}\right|^\varsigma d\mu_\theta =
\int g^2 \left(|\log \frac{g}{\int g^2d\lambda} + \beta d^\theta - \log Z|\right)^\varsigma d\lambda \nonumber \\
\leq  \int g^2 \left|\log \frac{g^2}{\int g^2d\lambda} \right|^\varsigma d\lambda 
+ \int g^2 \left(\beta  d^\theta\right)^\varsigma d\lambda
+ |\log Z|^\varsigma \int g^2  d\lambda \,\, \\
= \int g^2 \left|\log \frac{g^2}{\int g^2d\lambda} \right|^\varsigma d\lambda 
+ \beta^\varsigma \int f^2   d^{\theta\varsigma} d\mu_\theta
+ |\log Z|^\varsigma \int f^2  d\mu_\theta \nonumber 
\end{eqnarray}
Assume first that $\mu_\theta f^2 = \int g^2 d\lambda=1$. Then we have
\begin{eqnarray*}
\int g^2 \left(|\log \frac{g^2}{\int g^2d\lambda} |\right)^\varsigma d\lambda \leq
\int g^2 \left(\log_+ {g^2}\right)^\varsigma d\lambda  + D_\varsigma  \\
\leq  \left(\frac{2+\varepsilon}{\varepsilon}\right)^\varsigma 
\left(\log_+ \left(\int g^{2+\varepsilon}d\lambda\right)^{\frac2{2+\varepsilon}}\right)^\varsigma + D_\varsigma 
\end{eqnarray*}
with $D_\varsigma\equiv \sup_{x\in(0,1)}x\left(\log\frac1x\right)^\varsigma$. Choosing suitable $\varepsilon\in(0,1)$,
we can apply Sobolev inequality, (with constants $\bar C,\bar D\in(0,\infty)$), to get
\begin{eqnarray*}
\int g^2 \left|\log \frac{g^2}{\int g^2d\lambda} \right|^\varsigma \leq
\left(\frac{2+\varepsilon}{\varepsilon}\right)^\varsigma \left(\log_+ \left(\bar C \int |\nabla g|^2d\lambda + \bar D \int g^2 d\lambda \right)\right)^\varsigma   + D_\varsigma  \\
\leq  C_1\int |\nabla g|^2d\lambda + D_1
\end{eqnarray*}  
with 
$$
C_1\equiv s \left(\frac{2+\varepsilon}{\varepsilon}\right)^\varsigma \bar C 
$$
and
$$
D_1 \equiv  \left\{ s \left(\frac{2+\varepsilon}{\varepsilon}\right)^\varsigma \bar D + \gamma_{\varsigma,s} + D_\varsigma \right\}
$$
where $s\in(0,\infty)$ and  $\gamma_{\varsigma,s}\in(0,\infty)$ is a suitable constant.
Using the definition of $g$, we have
\begin{equation*}
\int |\nabla g|^2d\lambda \leq 2 \int |\nabla f|^2 d\mu_\theta + \frac12\beta^2\theta^2  \int f^2 d^{2(\theta-1)} d\mu_\theta
\end{equation*}
Now applying the U-bound of Theorem \ref{thm2.3}, we get

\begin{equation*}
\int |\nabla g|^2d\lambda \leq \left(2+ \frac12\beta^2\theta^2 C_\theta \right)\int |\nabla f|^2 d\mu_\theta + 
\frac12\beta^2\theta^2  D_\theta \int f^2 d\mu_\theta
\end{equation*}
Thus we get (for the normalised function $g$)
\begin{equation}\label{4.2.3}
\int g^2 \left|\log \frac{g^2}{\int g^2d\lambda} \right|^\varsigma \leq C_2 \int |\nabla f|^2 d\mu_\theta
+ D_2
\end{equation}
with some constants $C_2,D_2\in(0,\infty)$.
Now coming back to (\ref{4.2.1}), we note that since
$\theta\varsigma = 2(\theta-1)$, we can use again the U-bound of
Theorem \ref{thm2.3} to bound the second term from the right hand side of
this relation. Combining this with (\ref{4.2.3}), we arrive at
the following bound
 \begin{equation}\label{4.2.4}
\int f^2 \left|\log \frac{f^2}{\int f^2d\mu_\theta}\right|^\varsigma d\mu_\theta
\leq C \int |\nabla f|^2 d\mu_\theta + D
 \end{equation}
with the constants $C = C_2+ \beta^\varsigma C_\theta$ and $D= D_2+\beta^\varsigma D_\theta + |\log Z|^\varsigma$. 
At this stage we can remove the normalisation condition to arrive at the desired bound (\ref{4.2.a}).

\qed 
Using Theorem \ref{thm4.2}, we prove the following tight inequality.

\begin{theorem} \label{thm4.3}
For $\theta\in[1,2]$ and 
$\varsigma = \frac{2(\theta-1)}{\theta}$, let
$$
\Phi (x) \equiv x \left(\log\left(1+x\right)\right)^\varsigma 
$$
Under the assumption of Theorem \ref{thm4.2}, if additionaly $\mu_\theta$
satisfies Poincar\'e inequality, there exists a constant
$c_\theta\in(0,\infty)$ such that
\begin{equation} \label{4.3.a}
\mu_\theta \Phi(f^2) - \Phi(\mu_\theta f^2)  
\leq c_\theta  \int |\nabla f|^2 d\mu_\theta
\end{equation}

\end{theorem}

\textit{Proof}: 
First we note that
\begin{equation} \label{4.3.1}
\mu_\theta \Phi(f^2) - \Phi(\mu_\theta f^2) \leq  
\mu_\theta f^2 \left|\log \frac{1+f^2}{1+\mu_\theta f^2}\right|^\varsigma  
\end{equation}
and 
\begin{eqnarray} \label{4.3.2}
\mu_\theta f^2 \left|\log \frac{1+f^2}{1+\mu_\theta f^2}\right|^\varsigma  
= \mu_\theta \chi(f^2\geq \mu_\theta f^2)f^2 \left|\log \frac{1+f^2}{1+\mu_\theta f^2}\right|^\varsigma \\
+ \mu_\theta \chi(f^2\leq \mu_\theta f^2)f^2 \left|\log \frac{1+f^2}{1+\mu_\theta f^2}\right|^\varsigma \nonumber
\end{eqnarray}
On the set $\{f \geq \mu_\theta f^2\}$ we have $\frac{1+f^2}{1+\mu_\theta f^2}\leq \frac{f^2}{\mu_\theta f^2}$
and so
$$
\mu_\theta \chi(f^2\geq \mu_\theta f^2)f^2 \left|\log \frac{f^2}{\mu_\theta f^2}\right|^\varsigma \leq
\mu_\theta f^2 \left|\log \frac{f^2}{\mu_\theta f^2}\right|^\varsigma 
$$
On the other set $\{f \leq \mu_\theta f^2\}$, we have $\frac{1+\mu_\theta f^2}{1+ f^2}\leq 1+ \frac{\mu_\theta f^2}{ f^2}$,
and therefore
$$
\mu_\theta \chi(f^2\leq \mu_\theta f^2)f^2 \left|\log \frac{1+f^2}{1+\mu_\theta f^2}\right|^\varsigma 
\leq 2 \mu_\theta f^2
$$
Using these relations together with (\ref{4.3.2}) we have
\begin{equation} \label{4.3.3}
\mu_\theta \Phi(f^2) - \Phi(\mu_\theta f^2) \leq  
\mu_\theta f^2 \left|\log \frac{f^2}{\mu_\theta f^2}\right|^\varsigma + 2\mu_\theta f^2 
\end{equation}
and thus, by Theorem \ref{thm4.2}, we obtain
\begin{equation} \label{4.3.4}
\mu_\theta \Phi(f^2) - \Phi(\mu_\theta f^2) \leq  
C \mu_\theta |\nabla f|^2  + (D+2)\mu_\theta f^2 
\end{equation}
Now according to Lemma A.1 of \cite{LZ}, one has the following analog of Rothaus lemma for a probability measure
with Orlicz function $\Phi$ given in the theorem: 
$\exists a,b\in(0,\infty)\qquad $
\begin{equation} \label{4.3.5}
\nu\Phi(f^2) - \Phi(\nu f^2) \leq a\left[\nu\Phi((f-\nu f)^2) - \Phi(\nu (f-\nu f)^2)\right]
+ b \nu (f-\nu f)^2
\end{equation}
Combining (\ref{4.3.4}) and (\ref{4.3.5}) with the Poincar\'e inequality
for the measure $\mu_\theta$ 
$$
\mu_\theta(f-\mu_\theta f)^2 \leq \frac1M \mu_\theta|\nabla f|^2
$$
we arrive at the following result
$$
\mu_\theta \Phi(f^2) - \Phi(\mu_\theta f^2) \leq  
\left[aC + \frac{D+b}{M}\right]\, \mu_\theta |\nabla f|^2
$$

\qed

Summarising, in the current section in essence our methods were based on the fact that
the primary part of the interaction where a nice function of certain unbounded function
$d$ which length of the gradient $|\nabla d|$ (with respect to a given set of fields) was bounded from above and
stayed strictly away from zero. We also used number of times the Leibniz rule for the fields.

\subsection{From Coercive Inequalities to U-Bounds.}
\label{Sub4.3 LSq2U-Bounds}
For a probability measure $d\mu\equiv e^{-U}d\lambda/Z$, we have shown that
if for $q\in(1,2]$ the following bound is satisfied
$$
\int f^q \left(|\nabla U|^q + U\right) d\mu \leq C \int |\nabla f|^q d\mu + D \int |f|^q d\mu
$$ 
together with $q$-Poincar\'e inequality
$$
M \mu | f- \mu f|^q \leq \mu |\nabla f|^q \, ,
$$
then the following $\mathbf{LS}_q$ inequality holds
$$
\mu |f|^q \log \frac{|f|^q}{\mu |f|^q} \leq c \mu |\nabla f|^q
$$

We show that the following result in the converse direction is true as well.
\begin{theorem} \label{thm4.4}
Suppose $q\in(1,2]$ and for some constants $a,b\in(0,\infty)$ , we have
$$
|\nabla U|^q \leq aU +b
$$
and assume that the measure $d\mu\equiv e^{-U}d\lambda/Z$ satisfies
$\mathbf{LS}_q$. Then the following U-bound is true
$$
\int |f|^q U d\mu \leq C \int |\nabla f|^q d\mu + D \int |f|^q d\mu
$$
with some constant $C,D\in(0,\infty)$ independent of $f$.

\end{theorem}

\textit{Proof}: 
We note that by relative entropy inequality one has 
$$
\mu \left(|f|^q U \right) \leq \frac1\varepsilon \mu |f|^q\log \frac{|f|^q}{\mu |f|^q}
+ \left(\frac1\varepsilon \log \mu e^{\varepsilon U} \mu \right)\mu |f|^q 
$$
Hence, if $\mathbf{LS}_q$ is true, we get
$$
\mu \left(|f|^q U \right) \leq \frac{c}\varepsilon \int |\nabla f|^q d\mu
+ \left(\frac1\varepsilon \log \mu e^{\varepsilon U}  \right)\mu |f|^q 
$$
Thus we will be finished if we show $\mu e^{\varepsilon U}<\infty$. This follows from the following result.\\

\noindent \textbf{Exp-Bounds from $\mathbf{LS}_q$}\\
\begin{theorem} \label{thm4.5}
Assume that a measure $ \mu$ satisfies
$\mathbf{LS}_q$ with some $q\in(1,2]$.
Suppose that for some constants $a,b\in(0,\infty)$ , we have
$$
|\nabla f|^q \leq a f +b
$$
 Then the following exp-bound is true
$$
\mu e^{tf} < \infty
$$
 for all $t>0$ sufficiently small.
\end{theorem}
$\ $\\
\noindent{\textbf{Remark}}: For the case $q=2$ see \cite{AiSt}.\\

\textit{Proof}: 
By our assumption, we have
$$
\mu g^q\log\frac{g^q}{\mu g^q}\leq c \mu|\nabla g|^q  
$$
It is enough to prove the bound under additional assumption that
$f$ is bounded.  Namely,
given $L\in(0,\infty)$, replace $f$ by
$F \equiv \chi(|f|\leq L)f +L\chi(|f|> L)$. 
$F$ satisfies our assumptions with the same constants. So we will
get the claim letting $L$ go to $\infty$.

Since now $f$ is bounded, $\exp{tf}$ is integrable and we have
$$
\mu\left(e^{tf}\log\frac{e^{tf}}{\mu e^{tf}}\right) \leq
cq^{-q}t^q \mu \left(e^{tf} |\nabla f|^q \right) 
$$
By our assumption $|\nabla f|^q \leq a f +b$, so we get
$$
\mu\left(e^{tf}\log\frac{e^{tf}}{\mu e^{tf}}\right) \leq
caq^{-q}t^q \mu \left(e^{tf} f \right) + cbq^{-q}t^q \mu \left(e^{tf} \right)
$$
which can be rearranged to get
$$
(1-caq^{-q}t^{q-1})\mu\left(\frac{e^{tf}}{\mu e^{tf}}\log\frac{e^{tf}}{\mu e^{tf}}\right) \leq caq^{-q}t^{q-1} \log \mu \left(e^{tf} \right) + cbq^{-q}t^q
$$
Taking into the account that
$$
  \mu\left(\frac{e^{tf}}{\mu e^{tf}}\log\frac{e^{tf}}{\mu e^{tf}}\right) 
= t^2\frac{d}{dt}\frac1t\log  \mu e^{tf}   
$$
and setting $G(t)\equiv \frac1t\log  \mu e^{tf}$, after simple transformations 
we obtain the following differential inequality 
$$
\frac{d}{dt} G(t) \leq \beta t^{q-2 } G(t) +  \gamma t^{q-2 }
$$
with $\beta(t)\equiv  \frac{caq^{-q}}{(1-caq^{-q}t^{q-1})} $ and 
$\gamma (t)\equiv \frac{cbq^{-q}}{(1-caq^{-q}t^{q-1})} $
which are well defined for $caq^{-q}t^{q-1}<1$. Since $G(t)\to \mu f$ as $t\to 0$ and $q\in(1,2]$,
for $caq^{-q}t^{q-1}<\varepsilon<1$, 
after integration we get
$$
G(t)\leq \mu f + \frac{cbq^{-q}}{(q-1)(1-\varepsilon)} t^{q-1 }  + \frac{caq^{-q}}{(1-\varepsilon)}\int_0^t d\tau \ \tau^{q-2 } G(\tau) 
$$
In our range of $q\in(1,2]$, this can be solved by iteration.  Since $G(t)$ is nondecreasing, in this interval one also has
$$
G(t) \leq \mu f + \frac{cbq^{-q}}{(q-1)(1-\varepsilon)} t^{q-1 }  + \frac{caq^{-q}}{(q-1)(1-\varepsilon)} t^{q-1 } G(t) 
$$ 
which for $\frac{caq^{-q}}{(q-1)(1-\varepsilon)} t^{q-1 }\equiv \delta <1$ yields the following bound
$$
\mu e^{(1-\delta)tf}\leq \exp\{t\mu f + Ct^q\}
$$
with $C\equiv \frac{cbq^{-q}}{(q-1)(1-\varepsilon)}$.
One can check that our bound is independent of the cut off $L$ in the given interval of $t$. 
\\

By the above we have shown the equivalence of the $\mathbf{LS}_q$ and $U$-bounds
in particular in the cases of natural interactions dependent on the metric.
Similar considerations can be provided in the subquadratic case
for which the exponential bounds are known (see e.g. \cite{LO}, \cite{BCR}).

\section{Weighted U-Bounds and Coercive Inequalities.} 
\label{Sub4.4 WeightedU-Bounds}

Let $p\geq 2$ and suppose $f$ is a smooth function supported away from the origin. Starting with the identity
$$
d^{-\frac\alpha2}(\nabla f)e^{-\frac{\beta d^p}2} =  
d^{-\frac\alpha2}\nabla \left(fe^{-\frac{\beta d^p}2} \right)
+ \frac{p\beta}2 d^{p-\frac\alpha2-1} (\nabla d)fe^{-\frac{\beta d^p}2} \, ,
$$
squaring and integrating with the measure $d\lambda$, one obtains
\begin{eqnarray*}
\int d^{-\alpha}|\nabla f|^2 e^{-\beta d^p}d\lambda
\geq p\beta\int d^{p-\alpha-1} \nabla \left(fe^{-\frac{\beta d^p}2} \right)\cdot (\nabla d) fe^{-\frac{\beta d^p}2} d\lambda \nonumber \\
+  \frac{p^2\beta^2}4 \int d^{2p- \alpha-2} |\nabla d|^2 f^2 e^{-\beta d^p} d\lambda
\nonumber
\end{eqnarray*}
Hence, after integration by parts in the first term on the right hand side and simple rearrangements, one arrives at the following bound
\begin{eqnarray*}
\int d^{-\alpha}|\nabla f|^2 e^{-\beta d^p}d\lambda
\geq \frac{p^2\beta^2}4 \int  f^2 \left(d^{2p- \alpha-2} |\nabla d|^2 
\right)
e^{-\beta d^p} d\lambda \nonumber \\
- \int  f^2 \left[ \frac{p(p-\alpha -1)\beta}{2} d^{p-\alpha-2} |\nabla d|^2
+ \frac{p\beta}{2} d^{p-\alpha-1}  \Delta d \right] e^{-\beta d^p} d\lambda   \nonumber 
\end{eqnarray*}
If we choose $\alpha= p-2$ and assume $|\nabla d|\geq \frac1\sigma > 0$, we obtain
\begin{eqnarray*}
\int d^{-\alpha}|\nabla f|^2 e^{-\beta d^p}d\lambda
\geq \frac{p^2\beta^2}{4\sigma^2} \int  f^2  d^p  
e^{-\beta d^p} d\lambda \nonumber \\
-\int  f^2 \left[ \frac{p(p+1)\beta}{2}  |\nabla d|^2
+ \frac{p\beta}{2} d  \Delta d \right] e^{-\beta d^p} d\lambda   \nonumber 
\end{eqnarray*}
Finally assuming that there exists  constants $K\in(0,\infty)$ and $\delta\in(0,\frac{p^2\beta^2}{4\sigma^2})$, such that
$$
\frac{p(p+1)\beta}{2}  |\nabla d|^2
+ \frac{p\beta}{2} d  \Delta d \leq K + \delta d^p
$$
we arrive at 
$$
\left(\frac{p^2\beta^2}{4\sigma^2}-\delta\right) \int  f^2  d^p  
e^{-\beta d^p} d\lambda 
\leq \int d^{-\alpha}|\nabla f|^2 e^{-\beta d^p}d\lambda +  K \int f^2  e^{-\beta d^p}d\lambda 
$$
By adjusting the constant on the right hand side and replacing 
 $d^{-\alpha}$ by \par\noindent
$<d>^{-\alpha} \equiv (1+d^2)^{-\frac{\alpha}2}$, we conclude with the following result.

\begin{theorem}\label{thm4.4.1}
Let $d\mu\equiv e^{-\beta d^p}d\lambda /Z $ with $p > 2$. 
Suppose there are constants $\sigma\in[1,\infty)$ and  $K\in(0,\infty)$ and $\delta\in(0,\frac{p^2\beta^2}{4\sigma^2})$ such that
$|\nabla d| \geq \frac1\sigma$ and
$$
\frac{p\beta}{2}  |\nabla d|^2
+ \frac{p\beta}{2} d  \Delta d \leq K + \delta d^p
$$
Then there are constant $C,D\in(0,\infty)$ such that
$$
\mu  f^2  d^p  \leq 
C \mu  \left(<d>^{2-p}|\nabla f|^2 \right) 
+  D \mu  f^2  
$$
\end{theorem}
Using this bound, by similar arguments as in the proof of Poincar\'e inequality, (see Theorem \ref{thm3.1}), we now obtain 

\begin{theorem}\label{thm4.4.2}
Under the assumptions of Theorem \ref{thm4.4} there is a constant $M\in(0,\infty)$ such that
$$
M\ \mu \left(f-\mu f\right)^2 \leq    
\mu \left(<d>^{2-p}|\nabla f|^2 \right) 
$$
\end{theorem}
Finally following our strategy from the beginning of Section 4, (see proof of Theorem \ref{thm4.2}), with appropriate
amendments, we arrive at the following coercive inequality.

\begin{theorem}\label{thm4.4.3}
Under the assumptions of Theorem \ref{thm4.4} there is a constant $c\in(0,\infty)$ such that
$$
\mu \left(f^2 \log\frac{f^2}{\mu f^2}\right) \leq
c \mu \left(<d>^{2-p}|\nabla f|^2 \right) 
$$
\end{theorem}
 
\newpage

\subsection{Weighted U-Bounds and Coercive Inequalities: \\ Distributions with Slow Tails
on Riemannian Manifolds.} 
\label{Sub4.5 SlowTails_WeightedU-Bounds}

In this section we consider a noncompact smooth Riemannian manifold $\mathbb{M}$ of dimension $3\leq N<\infty$.  In this setup $d(x)$ denotes the Riemannian distance of a point $x$ from a given point $x_0\in\mathbb{M}$ called later on the origin. By $\nabla$ and $\Delta$ we denote the gradient and Laplace-Beltrami operators, respectively.
\\
The aim of this section is to discuss coercive inequalities involving probability measures $d\mu \equiv \rho dx$ with density 
(with respect to the corresponding Riemannian measure  $d\lambda$  on $\mathbb{M}$) which is 
 of the form
$\rho\equiv e^{-U(d)}/Z$ with leading part of the function $U$ given by a concave function (and therefore also defining a non-Riemannian distance on $\mathbb{M}$ ). In particular we will consider the following cases:\\ 
(i) 
$U(d)= \beta d^{\alpha}$, with $\alpha\in(0,\infty)$ and  $\beta>0$, \\
(ii) 
$U(d)=\beta \log(1 + d)$ with $\beta>0$.\\

Before we go on we recall the following Laplacian comparison theorem, (cf \cite{CLN}, \cite{Kas} (\cite{Q}, \cite{W1}-\cite{W2})).
For a complete Riemannian manifold $\mathbb{M}$ with $Ric\geq (N-1)K$ where $K\in\mathbb{R}$:\\
$(*)\phantom{*}$ If $K\leq 0$, then
$\Delta d\leq (N-1)d^{-1}+ (N-1)\sqrt{|K|}$\\ 
$(**)$ If $Ric \geq 0$, then
$\Delta d\leq (N-1)d^{-1}$\\

\noindent By similar computation as we have done in Section 2, for a smooth nonnegative function $f$ localised outside a ball $B_\varepsilon\equiv B_\varepsilon(x_0)$ centred at the origin we consider a field
\begin{equation} \label{4.5.0}
(\nabla f) e^{-U}= \nabla \left(f  e^{-U}\right)
 + f\left( U'\nabla d\right)e^{-U}.
\end{equation}
to which we will apply a functional 
 \begin{equation} \label{4.5.1}
\boldsymbol{\alpha}(\mathbf{v})\equiv \int W (\nabla d  \cdot \mathbf{v} )\  d\lambda.
\end{equation}
defined with a positive weight function $W\equiv W(d)$ to be specified later.
Using the fact that $|\nabla d |= 1$ (for $d\neq 0$), together with arguments involving H\"older inequality and integration by parts one arrives at the following bound
 
\begin{equation} \label{4.5.2}
\int f \mathcal{V}\  e^{-U}d\lambda \leq
\int W |\nabla f|e^{-U}d\lambda  
\end{equation}
with
$$
\mathcal{V} \equiv \chi_{\mathbb{M}\setminus B_\varepsilon}\left( W U' - 
\divergence(W\nabla d)\right)  
$$
Later on we will extend  $\mathcal{V}$ to $B_\varepsilon$ in a convenient way by adding an arbitrary bounded continuous function. 
\noindent One can handle a function of arbitrary sign replacing $f$ by
$|f|$ and using equality $\nabla |f| = \sgn(f)\nabla f$.
To include $f$ which are non-zero on a ball centred at the origin
we write $f = f_0 + f_1$ where $f_0 = \phi f$, $f_1 = (1-\phi)f$ and
$\phi(x) = \min(\varepsilon, \max(2\varepsilon - d(x),0))$.  Then
\begin{equation} \label{4.5.3}\begin{split}
\int   |f| \mathcal{V} d\mu  = \int_{d(x)\leq 2\varepsilon }|f| \mathcal{V} d\mu 
+ \int_{d(x)>2\varepsilon}|f| \mathcal{V} d\mu \\ 
\leq \sup_{\{d\leq 2\varepsilon\}}(\mathcal{V}) \int \phi |f|d\mu  + \int |f|_1 \mathcal{V} d\mu 
     \end{split}
\end{equation}
Next we have
\begin{equation} \label{4.5.4}
|\nabla f_1| \leq |\nabla f|+ \frac1\varepsilon \chi_{\{\varepsilon\leq d < 2\varepsilon\}}|f|,
\end{equation}
and therefore
\begin{equation} \label{4.5.5}\begin{split}
\int |f_1| \mathcal{V} d\mu 
\leq  \int W(1-\phi)|\nabla f | d\mu  +  
\sup_{\{\varepsilon\leq d < 2\varepsilon\}}( \varepsilon^{-1}W   ) 
\int_{\varepsilon\leq d < 2\varepsilon} |f| d\mu 
\end{split}
\end{equation}
Combining (\ref{4.5.0}) - (\ref{4.5.5}) we arrive at the following bound
\begin{equation} \label{4.5.6}
\int |f| \mathcal{V} d\mu 
\leq  \int W (1-\phi)|\nabla f | d\mu  +  
\sup_{\{d\leq 2\varepsilon\}}(\mathcal{V}) \int \phi |f| d\mu \phantom{AAAAAAAAAAAA}
\end{equation}
$$\phantom{AAAAAAAAAAAA}
+ \sup_{\{\varepsilon\leq d < 2\varepsilon\}}( \varepsilon^{-1} W ) \int_{\varepsilon\leq d < 2\varepsilon} |f| d\mu 
$$
Hence with 
$$B\equiv  \sup_{\{d\leq 2\varepsilon\}}(\mathcal{V}) +\sup_{\{\varepsilon\leq d < 2\varepsilon\}}( \varepsilon^{-1} W    )   ,
$$ 
we have
\begin{equation} \label{4.5.7}
\int |f| \mathcal{V} d\mu 
\leq  \int W  |\nabla f | d\mu  +  
B \int   |f| d\mu  
\end{equation}

\newpage

\textbf{Case (i)}\label{Case(i)}\\
\bigskip

\noindent For $U(d)= \beta d^{\alpha}$, with $\alpha\in(0,\infty)$ and  $\beta>0$,
choosing $W(d)= \alpha^{-1}d^\kappa$, with $\kappa\geq 1$, we have

\begin{equation}\label{4.5.i1}
\mathcal{V} \equiv W U' - \divergence(W\nabla d)=
U - \alpha^{-1}\kappa d^{\kappa-1} - \alpha^{-1} d^\kappa \Delta d  
\end{equation}
Thus if (*) holds, we have
\begin{equation}\label{4.5.i1a}
\mathcal{V}  \geq   \beta d^{\alpha-1+\kappa} - \chi_{\mathbb{M}\setminus B_\varepsilon} 
\left(   \alpha^{-1}\kappa N d^{\kappa-1}  
+ \alpha^{-1} (N-1)\sqrt{|K|}d^\kappa \right)
\end{equation}
Hence we conclude with the following result

\begin{theorem}\label{thm4.5.i1}
Let $d\mu\equiv e^{-U}d\lambda /Z $ with $U\equiv\beta d^\alpha$ where $\alpha\in (0,\infty)$. 
Suppose  $Ric\geq (N-1)K$ with $K\leq 0$.\\ 
$\bullet$ If $\alpha > 1$, then for any $\kappa \geq 1$, there exist constants $c_1,b_1\in(0,\infty)$
such that

\begin{equation} \label{4.5.i2}
\int |f| U d\mu 
\leq  c_1 \int  d^\kappa  |\nabla f | d\mu  +  
b_1  \int   |f| d\mu  
\end{equation}\\
$\bullet$ If $\alpha = 1$ and $\beta > \alpha^{-1}(N-1)\sqrt{|K|} $, then for any $\kappa \geq 1$, there exists constant $c_1,b_1\in(0,\infty)$
such that (\ref{4.5.i2}) is true.\\
$\bullet$ If $\alpha \in(0,1)$ and $Ric\geq 0$, then for any $\kappa \geq 1$, there exist constants $c_1,b_1\in(0,\infty)$
such that (\ref{4.5.i2}) is true.

Moreover if (\ref{4.5.i2}) holds, then for any $q\in(1,\infty)$, we have 
\begin{equation} \label{4.5.i3}
\int |f|^q U d\mu 
\leq  c_2 \int  d^{q(\kappa-\frac{\alpha}{p})}  |\nabla f |^q d\mu  +  
b_2  \int   |f|^q d\mu  
\end{equation}
with $c_2\equiv c_1\lambda q^{q-1}\beta^{\frac{q}{p}} [1- c_1/(p\lambda)]^{-1}$ and $b_2\equiv b_1[1- c_1/(p\lambda)]^{-1}$.
\end{theorem}
The second part follows from the first by substituting $f^q$ in place of $f$ and using
elementary arguments involving Young inequality.\\

As a consequence, by similar arguments as earlier in this section, we obtain the following result on possible coercive inequalities.
\begin{theorem}\label{thm4.5.1bis}
Let $d\mu\equiv e^{-U}d\lambda /Z $ with $U\equiv\beta d^\alpha$ where $\alpha\in (0,\infty)$. 
Suppose  $Ric\geq (N-1)K$ with $K\leq 0$.\\ 
$\bullet$ If $\alpha > 1$, then for any $\kappa \geq 1$, there exist constants $c\in(0,\infty)$
such that

\begin{equation}\label{4.5.i1b1b}
\mu |f|^q\log\frac{|f|^q}{\mu |f|^q} \leq
   c \int  d^{q(\kappa-\frac{\alpha}{p})}  |\nabla f |^q d\mu 
\end{equation}
\\
$\bullet$ If $\alpha = 1$ and $\beta > \alpha^{-1}(N-1)\sqrt{|K|} $, then for any $\kappa \geq 1$, there exist a constant $c \in(0,\infty)$
such that (\ref{4.5.i1b1b}) is true.\\
$\bullet$ If $\alpha \in(0,1)$ and $Ric\geq 0$, then for any $\kappa \geq 1$, there exist a constant  $c \in(0,\infty)$
such that (\ref{4.5.i1b1b}) is true.\\
As a consequence the following inequality holds
\begin{equation}\label{4.5.i1b1bis}
M\ \mu |f-\mu f|^q\leq
    \int  d^{q(\kappa-\frac{\alpha}{p})}  |\nabla f |^q d\mu 
\end{equation}
with some $M\in(0,\infty)$.
\end{theorem}

\textbf{Case (ii)}\label{Case(ii)}\\
\bigskip

\noindent
For $U(d)= \beta \log (1+d)$ with  $\beta>0$,
choosing $W(d)= d\log(1+d)$ and setting
\begin{equation}\label{4.5.ii1}
\begin{split}
\mathcal{V} \equiv U + \chi_{\mathbb{M}\setminus B_\varepsilon} \left(W \beta (1+d)^{-1} - \divergence(W\nabla d)\right) \\
=
U -  \chi_{\mathbb{M}\setminus B_\varepsilon}  [1+\log(1+d)]  - \chi_{\mathbb{M}\setminus B_\varepsilon}d\log(1+d) \Delta d  
\end{split}
\end{equation}
Thus if (*) holds, we have
\begin{equation}\label{4.5.ii1a}
\mathcal{V}  \geq U -   \chi_{\mathbb{M}\setminus B_\varepsilon} [1+\log(1+d)]  - \chi_{\mathbb{M}\setminus B_\varepsilon} d\log(1+d) \left[(N-1)d^{-1}+ (N-1)\sqrt{|K|}\right]
\end{equation}
Hence we conclude with the following result

\begin{theorem}\label{thm4.5.ii1}
Let $d\mu\equiv  (1+d)^{-\beta}d\lambda /Z $ with $\alpha\in (0,1)$. 
Suppose  $Ric\geq 0$. \\
If $\beta > N$, then

\begin{equation} \label{4.5.ii2}
\int |f| U d\mu 
\leq  c_1 \int  d\log(1+d)  |\nabla f | d\mu  +  
b_1  \int   |f| d\mu  
\end{equation}
with 
$$ c_1 \equiv\beta\cdot[\beta - N]^{-1}
$$
and
$$b_1\equiv \beta\cdot[\beta - N]^{-1}\cdot
\left(N + \sup_{\{d\leq 2\varepsilon\}}(\mathcal{V}) +\sup_{\{\varepsilon\leq d < 2\varepsilon\}}( \varepsilon^{-1} W    )  \right) ,
$$ 
Hence, there exist $c_q,b_q\in(0,\infty)$ such that
\begin{equation} \label{4.5.ii3}
\int |f|^q U d\mu 
\leq  c_q \int  d^q\log(1+d)  |\nabla f |^q d\mu  +  
b_q  \int   |f|^q d\mu  
\end{equation}
\end{theorem}
The second part follows from the first by substituting $f^q$ in place of $f$
and using the following Young inequality
$$
d  |\nabla f^q |= q\left(|f|^{q-1}\cdot d|\nabla f|\right)\leq 
\lambda^q d^q |\nabla f|^q +  \frac{q}p \lambda^{-p}|f|^q
$$
which implies
$$
\int  d\log(1+d)  |\nabla f^q | d\mu = \int d\log(1+d)q|f|^{q-1}|\nabla f| d\mu 
$$
$$
\leq
 \lambda^q\int d^q\log(1+d)  |\nabla f|^q d\mu 
+ \frac{q}p\lambda^{-p}\int  \log(1+d)q|f|^{q }  d\mu
$$
From this and (\ref{4.5.ii2}), choosing $c_1\frac{q}{p}\lambda^{-p}<1$, one obtains
$$
\int |f|^q U d\mu 
\leq  c_q \int  d^q\log(1+d)  |\nabla f |^q d\mu  +  
b_q  \int   |f|^q d\mu  
$$
with 
$c_q\equiv c_1\lambda^q(1-c_1\frac{q}{p}\lambda^{-p})^{-1}$ and $b_q\equiv b_1(1-c_1\frac{q}{p}\lambda^{-p})^{-1}$.


\noindent As a consequence of the above theorem, using arguments similar to those of sections 4.1 and  4.2, we derive the following result on possible coercive inequalities.
\begin{theorem}\label{thm4.5.2c}
Let $d\mu \equiv e^{-\beta \log(1+d)}dx/Z $ with $\beta >N $. Suppose $Ric\geq 0$.
Then for any $q\geq 1$, there are constants $M_q,c_q\in(0,\infty)$, such that
\begin{equation}\label{4.5.ii1_c}
M_q \mu |f-\mu f|^q\leq \mu (1+d)^q\log(e+d)|\nabla f|^q
\end{equation}
and 
\begin{equation}\label{4.5.ii2_c}
 \mu |f|^q\log \frac{|f|^q}{\mu |f|^q}\leq c_q \mu (1+d)^q\log(e+d)|\nabla f|^q
\end{equation}
\end{theorem}
\noindent \subsubsection{\textit{Weighted Inequalities at Large $\beta$}}  \label{Weighted Ineqs_ii}

Let $U\equiv \beta \log(1+d)$, with $\beta>N\equiv dim(\mathbb{M})$. 
While the above results are true for any $\beta>N$, we will show that for sufficiently
big $\beta$ and $Ric\geq 0$ due to the special nature of the interaction it is possible to improve the weight in the Poincar\'e and related Log-Sobolev inequalities.\\
We start from noting that for a nonnegative differentiable function supported outside a ball
of radius $r$  centred at the origin, one has 
$$
\int (1+d)|\nabla f| e^{-U}dx \geq \int (1+d)\nabla d\cdot \nabla f e^{-U}dx
$$
$$
= \int (1+d)\left[\nabla d\cdot \nabla\left( f e^{-U} \right)+ f\nabla d\cdot \nabla U\right]dx 
$$
and so, taking into the account that $|\nabla f|^2=1$, one gets
$$
\int f \left[ \beta  - 1 - (1+d) \Delta d  \right]  e^{-U}  dx \leq \int (1+d)|\nabla f| e^{-U}dx
$$
When $Ric\geq 0$, we have  $\Delta d\leq (N-1)d^{-1}$ which implies the following bound
\begin{equation}\label{4.5.1.1}
M_\beta  \int f  e^{-U}  dx \leq \int (1+d)|\nabla f| e^{-U}dx
\end{equation}
where $M_\beta\equiv\left[ \beta  - N -\frac{(N-1)}r) \right]$
Since $|\nabla f|\geq |\nabla |f| |$, this inequality remains true for
not necessarily positive function with $f$ replaced by $|f|$ on the right
hand side.
Let now consider the following cutoff function
$$ 
\chi(t)\equiv\begin{cases}
              1\ for\ 0\leq t\leq 2r\\
              1-\frac{(t-r)}{L} )\ for\ 2r\leq t\leq R\\
              0\ for \ t\geq R  
             \end{cases}
$$
with some $R>2r$ to be chosen later.
Setting $\tilde f_1 \equiv (f-\mu f)\chi$ and $\tilde f_2 \equiv (f-\mu f)\chi$, we have
$$
\mu |f-\mu f|\leq \mu |\tilde f_1| + \mu |\tilde f_2| 
$$
As $\tilde f_1$ is compactly supported Lipschitz function, there is an $m\equiv m_R\in(0,\infty)$
independent of the function $f$, such that
$$
\mu |\tilde f_1| \leq m_R^{-1}\mu |\nabla \tilde f_1|\leq  m_R^{-1}\mu(|\nabla f|\chi) 
+  \frac{1}{m_R(R-2r)}\mu\left(|f-\mu f|\chi({2r<d<R})\right)
$$
The second term on the right hand side can be treated with the help of (\ref{4.5.1.1}) as follows.
Setting $\hat \chi$ to be a Lipschitz extension of $\chi({2r<d<R})$ supported outside the ball of radius $r$, we have 
$$\mu\left(|f-\mu f|\chi({2r<d<R})\right)\leq \mu\left(|f-\mu f|\hat\chi \right)
\leq M_\beta^{-1}\  \mu (1+d)|\nabla f| 
$$
$$+ M_\beta^{-1} \sup|\nabla\hat\chi|\  \mu |f-\mu f|
$$
Thus we obtain
\begin{equation}\label{4.5.1.2}
 \mu |\tilde f_1| \leq m_R^{-1}\mu |\nabla \tilde f_1|\leq   
\left[m_R^{-1}
+  \frac{1}{m_R(R-2r)}M_\beta^{-1}\right] \mu (1+d)|\nabla f| \hat\chi
\end{equation}
$$+
\frac{1}{m_R(R-2r)} M_\beta^{-1} \sup|\nabla\hat\chi| \mu |f-\mu f| 
$$
On the other hand applying (\ref{4.5.1.1}) to $\tilde f_2$ we obtain
\begin{equation}\label{4.5.1.3}
\mu |\tilde f_2| \leq M_\beta^{-1} \mu (1+d)|\nabla f|(1-\chi) +  M_\beta^{-1}\frac{1+R}{R-r}\mu\left(|f-\mu f|\chi({r<d<R}) ) \right)
\end{equation}
Combining (\ref{4.5.1.2}) and (\ref{4.5.1.2}) we arrive at
\begin{equation}\label{4.5.1.4}
\mu |f-\mu f|\leq  c_0 \mu (1+d)|\nabla f| 
+
b_0 \mu |f-\mu f|  
\end{equation}
with
$$
c_0\equiv\left[m_R^{-1}+  (\frac{1}{m_R(R-2r)}+1)M_\beta^{-1}\right]
$$
and
$$
b_0\equiv  M_\beta^{-1}\left(\frac{1}{m_R(R-2r)}\sup|\nabla\hat\chi|+ \frac{1+R}{R-r}\right) 
$$
Since given $R>2r$, one can choose $\beta >N$ sufficiently large so that $b_0<1$, we conclude with the following result

\begin{theorem}\label{thm4.5.1.1}
Suppose $U=\beta \log(1+d)$, with $\beta>N$, and $Ric\geq 0$. Then there exists $\beta_0 > N$, such that for any $\beta >\beta_0$, one has 
\begin{equation}\label{4.5.1.5}
M \mu |f-\mu f|   \leq \mu (1+d)|\nabla f| 
\end{equation}
with some constant $M\in(0,\infty)$ independent of $f$.
Consequently, we have
\begin{equation}\label{4.5.1.5a}
M_q \mu |f-\mu f|^q   \leq \mu (1+d)^q|\nabla f|^q 
\end{equation}
with some constant $M_q\in(0,\infty)$
\end{theorem}
The second part of the theorem follows by similar arguments as the ones used in the proof
of Proposition 2.3 in \cite{BZ}.\\
Next we study the relative entropy estimate as follows. For a non-negative function $f$, setting 
$f_1 \equiv f \chi$ and $f_2\equiv f(1-\chi)$
with the same Lipschitz cutoff  function $\chi$ , we have
$$
\mu f\log \frac{f}{\mu f} \leq \mu f_1\log \frac{f_1}{\mu f_1} + \mu f_2\log \frac{f_2}{\mu f_2}
$$
Since the function $f_1$ is compactly supported and the density of the measure $\mu$ restricted to the ball $B_R(x_0)$ bounded and bounded away from zero, (via the arguments involving Sobolev inequality) we get
\begin{equation}\label{4.5.1.6}
\mu f_1\log \frac{f_1}{\mu f_1} \leq c_1 \mu |\nabla f_1| \leq c_1 \mu (|\nabla f |\chi)
+ b_1\sup|\nabla \chi| \mu f
\end{equation}
with some constants $c_1,b_1\in(0,\infty)$ independent of $f$.
Next we apply similar arguments based on Sobolev inequality with the function $F\equiv \frac{f_2e^{U}}{\int f_2e^{U}dx}$ and the Riemannian measure $dx$ to get
$$
\int F\log \frac{F}{\int F dx} dx \leq a \int |\nabla F|dx + b \int F dx 
$$
with some constants $a,b\in(0,\infty)$.
Hence we have 
\begin{equation}\label{4.5.1.7}
\mu f_2\log\frac{f_2}{\mu f_2}  \leq a \mu |\nabla f| (1-\chi) +  \mu f(1-\chi) (a|\nabla U| + b-\log Z)
+ \mu f_2 U 
\end{equation}
In our current setup we have $|\nabla U|\leq \beta$. Moreover, by simple relative entropy arguments, we have
$$
\mu f_2 U 
= \frac1\lambda \mu f_2 \log \frac{e^{\lambda U}}{\mu e^{\lambda U}}  
+ \frac1\lambda\log \mu e^{\lambda U} \mu f_2 
$$
$$\leq \frac1\lambda \mu f_2 \log\frac{f_2}{\mu f_2} 
+ \frac1\lambda\log \mu e^{\lambda U} \mu f_2 
$$
which hold provided that $\beta>N+\lambda$. If we can choose $\lambda > 1$, this together with (\ref{4.5.1.7}) implies
\begin{equation}\label{4.5.1.8}
\mu f_2\log\frac{f_2}{\mu f_2}  \leq c_2\mu |\nabla f| (1-\chi) +  b_2\mu f(1-\chi) 
\end{equation}
with
$$
c_2 \equiv a(1-\lambda^{-1})^{-1} 
$$
and
$$
b_2 \equiv (1-\lambda^{-1})^{-1} \left[a\beta + b-\log Z
 \frac1\lambda\log \mu e^{\lambda U} \right]
$$
Combining (\ref{4.5.1.8}) and (\ref{4.5.1.6}) we arrive at the following result

\begin{theorem}\label{thm4.5.1.2}
Suppose $U=\beta \log(1+d)$, with $\beta>N$, and $Ric\geq 0$. Then there exists $\beta_0 > N$, such that for any $\beta >\beta_0$, one has 
\begin{equation}\label{4.5.1.9}
\mu f\log\frac{f}{\mu f}  \leq \bar c \mu (1+d) |\nabla f| + \bar b \mu f
\end{equation}
with some constant $\bar c, \bar b \in(0,\infty)$ independent of $f$.
Consequently, if the weighted Poincar\'e inequality (\ref{4.5.1.5a}) is true for $q>1$, we have
\begin{equation}\label{4.5.1.9a}
\mu f^q\log\frac{f^q}{\mu f^q}  \leq  c_q \mu (1+d)^q |\nabla f|^q 
\phantom{AAAAAAAAAAAAAAAAAA}(\textbf{\textrm{WLS}}_q) \nonumber
\end{equation}
with some constant $c_q \in(0,\infty)$. 
\end{theorem}
We remark that (\ref{4.5.1.9}) implies similar weighted $LS_q$ inequality with $f$ replaced by $|f|^q$ and $|\nabla f|$ by its $q$-th power (which follows simply by substitution and use of H\"older inequality), while the tightening is obtained via Rothaus arguments (see e.g. \cite{BZ}).

\section{Optimal control distance on the Heisenberg Group. }
\label{5.Optimal control distance} 

Heisenberg group $H_l$ as a manifold is isomorphic to 
${\cal R}^{2l+1} = {\cal R}^{2l}\times {\cal R}$ with the
multiplication given by the formula
$$(x_1, z_1)\circ(x_2, z_2) = 
(x_1 + x_2, z_1 + z_2 + {\frac 1 2}S(x_1, x_2) )$$
where $S(x, y)$ is standard symplectic form on ${\cal R}^{2l}$:
$$S(x, y) = \sum_{i=1}^l (x_iy_{i+l} - x_{i+l}y_i).$$
Vector fields spanning the corresponding Lie algebra are give as follows
$$X_i = \partial_{x_i} + {\frac 1 2}x_{i+l}\partial_z,$$
$$X_{i+l} = \partial_{x_{i+l}} - {\frac 1 2}x_i\partial_z,$$
$$Z = \partial_z$$
where $i = 1,\dots, l$.

More generally, we say  that a
Lie algebra ${\bf n}$ is a stratified Lie algebra if it can be written as
$$
{\bf n} = \oplus_{i}^{m}{\bf n}_i,$$
$$[{\bf n}_i, {\bf n}_j] \subset {\bf n}_{i+j}$$
and ${\bf n}$ is generated by ${\bf n}_1$.  
Note that stratified Lie algebra is nilpotent.

We say that Lie group $N$ is stratified if it is connected, simply
connected and its Lie algebra ${\bf n}$ is stratified.  Since for stratified
groups exponential mapping is a diffeomorphizm from ${\bf n}$ to $N$,
one can identify $N$ with ${\bf n}$.

A Lie algebra is step two if it is stratified with $m = 2$.  In other words
it can be written in the form
$$
{\bf n} = {\bf v} \oplus {\bf z}$$
where ${\bf z}$ is the center (that is $[{\bf n}, {\bf z}] = {0}$)
and $[{\bf v}, {\bf v}] \subset {\bf z}$.

On a stratified Lie algebra ${\bf n}$ we define dilations by the formula
$$\delta(s) x = s^ix$$
for $x \in {\bf n}_i$ (and extend linearly to the whole ${\bf n}$.
For $s \ne 0$ $\delta(s)$ is an automorphism of ${\bf n}$.  One
can also define dilations on the corresponding group: 
$\delta(\exp(X)) = \exp(\delta(X))$.

A Lie algebra ${\bf n}$ is of H-type (Heisenberg type) if it is step two
and there exists an inner product $\langle \cdot, \cdot \rangle$ on
${\bf n}$ such that ${\bf z}$ is an orthogonal complement to
${\bf v}$, and the map $J_Z:{\bf v} \mapsto {\bf v}$ given by
$$
\langle J_Z X, Y\rangle = \langle [X, Y], Z\rangle
$$
for $X, Y \in {\bf v}$ and $Z \in {\bf z}$ satisfies $J_Z^2 = -|Z|^2I$
for each $Z \in {\bf z}$.  Equivalently,
for each $v\in {\bf v}$ of length $1$ the mapping $ad_{v}^{*}$
given by
$$
\langle ad_{v}^{*}z, y\rangle = \langle z, ad_v y \rangle
= \langle z, [v, y] \rangle 
$$
is an isometry from ${\bf z}^{*}$ into ${\bf v}^{*}$.

An H-type group is a connected and simply connected Lie group $N$
whose Lie algebra is of H-type.  We can identify H-type group $N$ with
its Lie algebra ${\bf n}$ defining multiplication on ${\bf n}$ by
the formula:
$$
(v_1, z_1)\cdot(v_2, z_2) = (v_1+v_2, z_1+z_2 + \frac{1}{2}[v_1, v_2])$$
where $v_1, v_2 \in {\bf v}$ and $z_1, z_2 \in {\bf z}$.

It is easy to see that Heisenberg group is an H-type group.  Also
H-type group with one-dimensional center is isomorphic to
the Heisenberg group, however
there exist H-type groups with center of arbitrary high dimension
\cite{Ka}.

On H-type group we consider vector fields $X_1,\dots,X_n$ which form
an orthonormal basis of ${\bf v}$ and we introduce the following operators

Subelliptic gradient:
$$\nabla f = (X_1f, \dots, X_{n}f)$$

Kohn laplacian
$$
\Delta = \sum_{i=1}^{n} X_i^2.$$
On Heisenberg group $H_l$ $n=2l$ and
$$
\Delta = \sum_{i=1}^{2l}\partial_{x_i}^2 
+\partial_z\sum_{i=1}^{l}(x_{i+l}\partial_{x_i}-x_i\partial_{x_{i+l}})+
\frac{|x|^2}{4}\partial_z^2.
$$
On general H-type group we similar, but more complicated expression:
$$
\Delta = \sum_{i=1}^{n}\partial_{v_i}^2 
+ \sum_{i=1}^{k}\partial_{z_i}\sum J_{\alpha,i}
+ \frac{|v|^2}{4}\sum_{i=1}^{k}\partial_{z_i}^2
$$
where $J_{\alpha, i}$ are vector fields corresponding to rotations.

Length of a curve: smooth $\gamma: [0, 1] \mapsto G$ is admissible
if 
$\gamma'(s) = \sum _{i=1}^{n} a_i(s)X_i(\gamma(s))$.
If $\gamma$ is admissible,
then $|\gamma| = \int_0^1 (\sum_{i=1}^{n} a_i^2(s))^{1/2}$.

Distance
$$d(g) = \inf |\gamma|$$
where infimum is taken over all admissible $\gamma$ such that $\gamma(0) = e$ and
$\gamma(1) = g$.

$d$ is homogeneous of degree $1$ with respect to the dilations $\delta(s)$,
namely for $s>0$
$$d(\delta(s)g) = sd(g).$$

\begin{lemma}\label{lem6.1}
On H-type group $Z$ distance $d((v,z))$ depends only on $|v|$ and $|z|$.
Moreover
if $\bar v, \bar z \in H_1$, $|v| = |\bar v|$, $|z| = |\bar z|$, then
$d((v, z)) = d((\bar v, \bar z))$.
\end{lemma}

\textit{Proof}:  Fix vectors $V, Z \in N$ such that $|V| = 1$, $|Z| = 1$,
$v = |v|V$, $z = |z|Z$.  Put $X = J_Z(V)$.  Since $J_Z$ is antisymmetric
and $J_Z^2 = I$, $J_Z$ is orthogonal, so $|X| = 1$.  Also, for
any $S \in {\bf z}$ of length $1$, we have
$$|\langle [X, Y], S \rangle| = |\langle J_SX, Y\rangle|
\leq |X||Y|$$
so since 
$$
\langle [V, X], Z \rangle = \langle J_Z V, X \rangle =
\langle X, X \rangle = |X|^2 = 1$$
we have $[V, X] = Z$.

Now, it is easy to see that the subgroup (in fact a subspace)
of $N$ generated by $V, X, Z$ is isomorphic to $H_1$.  Consequently,
using images of curves from $H_1$ to join $)$ with $(v, z)$ we
see that $d((v,z)) \leq d(((|v|, 0), z))$ where on the right hand
we have distance in $H_1$.

To get inequality in the opposite direction consider quotient group
$N/M$ where $M = \{t \in {\bf z}: \langle t, Z\rangle = 0\}$.  It
is easy to see that $N/M$ is still an H-type group (note that
since $N/M$ has one dimensional center it is enough to check the
defining property just for $J_Z$).  Hence, $N/M$ is isomorphic
to the Heisenberg group of appropriate dimension.  For
Heisenberg group our claim is well-known.
\qed

If is known \cite{Monti} that on Heisenberg group
if $g = (x, z)$ and $x \ne 0$ then
$d$ is smooth at $g$ and $|\nabla d| = 1$, however when $x = 0$
than $d$ is not differentiable at $g$.

\begin{lemma}
Let $A_\epsilon = {(r, z) \in {\cal R}^2: z >0, r > -\epsilon z}$. 
There is $\epsilon>0$ and a smooth function $\psi(r, z)$ defined on
$A_\epsilon$ such that on each group $N$ of H-type
$$d((x,z)) = \psi(|x|, |z|).$$
Moreover, $\partial_r \psi < 0$ when $r = 0$.
\end{lemma}

\textit{Proof}: First, by Lemma \ref{lem6.1} without loss of generality
we may assume that $N= H_1$.  Also, if $|x_1| = |x_2|$ 
and $|z_1| = |z_2$, then
$d(x_1, z_1) = d(x_2, z_2)$, so $\psi$ is uniquely defined
for $r\geq 0$.  We need to show that it has smooth extension
to $A_\epsilon$.  Since $d$ is homogeneous, it is enough to
construct smooth extension in a neighbourhood of a single
point $g = (0, 1)$.  

There exist a smooth
geodesic (length minimizing curve)
$\gamma$ joining $e = (0,0)$ and $g$.  We use length as
a parametrization of $\gamma$, so $\gamma(d(g)) = g$.  For $s < s_0 = d(g)$
we have $d(\gamma(s)) = s$.  

Let $\gamma(s) = (\gamma_x(s), \gamma_z(s))$.  Since square of Euclidean
distance is smooth $|\gamma_x|^2$ is smooth.  We can write
$|\gamma_x|^2(s) = (s-s_0)^2\rho(s)$ where $\rho$ is smooth and
$\rho(s_0) = 1$, so $|\gamma_x|^2(s)$ has a square root
$\phi(s) = (s_0-s)\rho^{1/2}(s)$ which is smooth for $s$ close to $s_0$.
Since both $\phi$ and $|\gamma_x|$ are positive square roots of
$|\gamma_x|^2$ for $s_0 - \epsilon < s<s_0$ we have
$$
|\gamma_x(s)| = \phi(s)$$
for $s_0 - \epsilon < s \leq s_0.$
Put
$$
\eta(s, t) = (t\phi(s), t^2\gamma_z(s)).$$ 
Since $\gamma$ is admissible
$|\gamma_z|'(s_0) = 0$ so the Jacobi matrix at $(s, t) = (s_0, 1)$ is
$$
\left(
    \begin{array}{cc}
     -1 & 0 \\
     0 & 2 \\
     \end{array}
 \right)
$$
and by the inverse function theorem $\eta$ is invertible in
a neighbourhood of $(s_0, 1)$.  So, there exist $f_1, f_2$ such that
$$(r, p) = \eta(f_1(r,p), f_2(r,p)).$$
We claim that $\psi(r,p) = f_1(r,p)f_2(r,p)$ give us extension of
$\psi$ to a neighbourhood of $g$.  Consider $(x, z)$ close to
$g$.  Let $(s, t) = (f_1(|x|, z), f_2(|x|, z))$.  We have
$$
|x| = t\phi(s) = t|\gamma_x(s)| = |(\delta_t\gamma(s))_x|,$$
$$
z = t^2\gamma_z(s) = (\delta_t\gamma(s))_z$$
so
$$
d((x, z)) = d(\delta_t\gamma(s)) = td(\gamma(s)) = ts = 
f_1(r,z)f_2(r,z) = \psi(r,z).$$

Now it remains to find sign $(\partial_r\psi)(0, z)$.  Form
equality $(r, p) = \eta(f_1(r,p), f_2(r,p))$ we see
$I = \eta'\cdot f'$.  We substitute $(r, p) = (0, 1)$ and
note that this corresponds to $(s_0, 1)$.  So
$$
\left(
    \begin{array}{cc}
     1 & 0 \\
     0 & 1 \\
     \end{array}
 \right)
=
\left(
    \begin{array}{cc}
     -1 & 0 \\
     0 & 2 \\
     \end{array}
 \right)\cdot
\left(
\begin{array}{c}
\partial_r f \\
\partial_p f \\
\end{array}
\right)
$$
and using first row we get
$ 1 = -(\partial_r f_1)(0, 1)$, $0 = -(\partial_r f_2)(0, 1)$
so
$$(\partial_r \psi)(0,1) = (\partial_r f_1)(0,1)f_2(0,1) 
+ f_1(0,1)(\partial_r f_2)(0,1) = (\partial_r f_1)(0,1) = -1$$ 
\qed

\begin{theorem}\label{thmhd}
If $N$ is an H-type group, then
there is $K$ such that if $d(g) \geq 1$, then
$$\Delta d \leq K$$
where $\Delta$ is understood in the sense of distributions.
\end{theorem}

\textit{Proof}:  Due to homogeneity, it is enough to prove the
inequality only for $g$ with $d(g) = 1$ (more precisely, in a small
neighbourhood of each such $g$).  Namely, if $s = d(g) > 1$
then
$$
\Delta d(g) = s^{-2}\Delta d(\delta(s)g) = s^{-1}\Delta d(g).$$
Next, $d((x, z)$ is smooth when $x \ne 0$, so it is enough to
prove the inequality in a small neighbourhood of $(0,z_0)$ where
$z_0>$ is chosen so that $d((0, z_0)) = 1$.  

Below we give computation on Heisenberg group:
$$
\partial_{x_i} d((x, z)) = \partial{x_i}\psi(|x|, z) 
= \frac{x_i}{|x|}\partial_r\psi(|x|, z),$$
$$
\partial_{x_i}^2 d((x, z)) = \partial{x_i}(\frac{x_i}{|x|}\partial_r\psi(|x|, z))$$
$$
= \frac{x_i^2}{|x|^2}\partial_r^2\psi(|x|, z)+ (\frac{1}{|x|}\partial_r\psi(|x|, z)-\frac{x_i^2}{|x|^3}\partial_r\psi(|x|, z),$$
$$
\sum_{i=1}^{2n}\partial_{x_i}^2 d((x, z)) = \frac{2n-1}{|x|}\partial_r\psi(|x|, z)+\partial_r^2\psi(|x|, z),$$
$$
(x_{i+n}\partial_{x_i}-x_i\partial_{x_{i+n}})d((x, z)) = 
(\frac{x_{i+n}x_i}{|x|} - \frac{x_ix_{i+n}}{|x|})\partial_r\psi(|x|, z) = 0,$$
$$
\Delta d((x,z)) =  \frac{2n-1}{|x|}\partial_r\psi(|x|, z)+\partial_r^2\psi(|x|, z) + \frac{|x|^2}{4}\partial_z^2\psi(|x|, z).$$
Since $\psi$ is smooth the second term and third term is bounded in
a neighbourhood
of $(0, z_0)$.  Since $\partial_r\psi(0, z_0) < 0$ the first term is
unbounded, but negative in a neighbourhood of $(0, z_0)$, which gives
the claim on Heisenberg group. 

On general H-type groups instead of
$x_{i+n}\partial_{x_i}-x_i\partial_{x_{i+n}}$ one must handle
the $J_{\alpha,i}$ term.  However, since $J_{\alpha,i}$ generates
rotations in $v$ space and $d$ is rotationally invariant again
$J_{\alpha,i}d = 0$.
\qed

\subsection{Counterexample for homogeneous norm}

On stratified groups $N$ one may introduce a homogeneous norm, that is
a continuous function $\phi:N \mapsto [0,\infty)$ such that $\phi(e) = 0$,
$\phi(x) > 0$ for $x\ne e$ and
$\phi(\delta_s(x)) = s\phi(x)$ for $s>0$.  Homogeneous norms are
equivalent to each other, if $\phi_1$ and $\phi_2$ are two
homogeneous norms, then there is $C$ such that
$$
C^{-1}\phi_1 \leq \phi_2 \leq C\phi_1.$$

The optimal control distance $d$ gives one example of homogeneous norm,
but there are others.  In particular, it is possible to choose
homogeneous norm so that it is smooth for $x \ne e$ (we will call
such homogeneous norm {\it smooth}).  Smooth homogeneous norms are
convenient in many situations.  For smooth  homogeneous norm $\phi$
the condition $(\Delta \phi)(x) \leq K$ for $\phi(x) \geq 1$ is
automatically satisfied.  However, we are going to prove that 
for such norm $|\nabla \phi|(x) = 0$ for some $x\ne e$, and consequently
log-Sobolev inequality like the one for optimal control distance can not
hold.

\begin{theorem}
Let $N$ be a stratified group, and $\phi$ be a smooth homogeneous
norm on $N$.  There exists $x\ne e$ such that $|\nabla \phi|(x) = 0$.
\end{theorem}

{\it Proof}:  Let $X_1,\dots,X_n$ be a basis of ${\bf n}_1$.  We
claim that for $(a_1,\dots,a_n) \in R^n -\{0\}$,
\begin{equation} \label{homodiff}
\sum a_i(X_i\phi)(\exp(\sum a_iX_i)) > 0.
\end{equation}
Namely, $\exp(t\sum a_iX_i)$ is a one parameter subgroup of $N$, so
$$
\partial_t(\phi(\exp(t\sum a_iX_i)) = \sum a_i(X_i\phi)(\exp(t\sum a_iX_i))$$
However, by homogeneity
$$
\partial_t(\phi(\exp(t\sum a_iX_i)) = \partial_t(t\phi(\exp(\sum a_iX_i))
= \phi(\exp(\sum a_iX_i) > 0$$
so (\ref{homodiff}) holds.

Using the $X_1,\dots,X_n$ basis we
identify ${\bf n}_1$ with $R^n$.  This identification gives us scalar
product on ${\bf n}_1$.  We extend this scalar product to a
scalar product on ${\bf n}$ such
that ${\bf n}_i$ is orthogonal to ${\bf n}_j$ for $i \ne j$.

Let $S$ ($\tilde S$) be the unit sphere
in ${\bf n}_1$ (in ${\bf n}$ respectively).  
Define mapping $\eta: S \mapsto S$ by the formula
$\eta(x) = \frac{(\nabla\phi)(\exp(x))}{|\nabla \phi|(\exp(x))}$
(note that we use identification ${\bf n}_1 = R^n$ here).
By (\ref{homodiff}) on $S$ $|\nabla \phi|(\exp(x)) > 0$ so $\eta$ is well
defined.  Also, $\eta$ is homotopic with identity.  Namely put
$\chi(\sum a_iX_i) = (a_1,\dots,a_n)$.
If $f_t$ is defined by the formula $f_t(x) = t\eta(x) + (1-t)\chi$, 
then for $x = \sum a_iX_i$ we have $\langle f_t(x), x\rangle > 0$, so
$f_t$ takes values in $R^n-\{0\}$.  Consequently
$g_t(x) = \frac{f_t(x)}{|f_t(x)|}$ gives homotopy of mappings
from $S$ to $S$.

If $(\nabla\phi)(\exp(x)) \ne 0$ on $\tilde S$,
then $\eta$ is homothopic to a constant.  Namely, $\tilde S$ contains
a homeomorphic copy of $n+1$ dimensional disc $D$ having $S$ as a boundary
and $\frac{(\nabla\phi)\circ\exp}{|(\nabla\phi)\circ\exp|}$ gives
required homotopy.  However, it is well known that identity of
the sphere is not homotopic to a constant -- so we reach
contradiction with assumption that $(\nabla\phi)(\exp(x)) \ne 0$.
\qed

\begin{lemma}\label{hom-taylor}
If $f$ is smooth function on a stratified group $N$, $d$ is optimal
control metric on $N$, $x_0\in N$ is fixed then
$$|f(x)-f(x_0)| \leq O(d(x,x_0)).$$
If additionally $(\nabla f)(x_0) = 0$, then
$$|f(x)-f(x_0)| \leq O(d^2(x,x_0)).$$
\end{lemma}

{\it Proof}: Let $\gamma: [0,1] \mapsto N$ be an admissible curve joining
$x_0$ and $x$. We have $\gamma'(s) = \sum a_i(s)X_i(\gamma(s))$, so
$$
|f(x)-f(x_0)| = \int_0^1 |(f\circ\gamma)'| = 
\int_0^1 |\sum a_i(s)(X_if)\circ\gamma|$$
$$
\leq  \int_0^1 |\gamma'||(\nabla f)\circ \gamma| 
\leq |\gamma|\sup_{s\in[0,1]}|(\nabla f)\circ \gamma(s)|.$$
Put $r = d(x,x_0)$. If $|\gamma|\leq r+\varepsilon$, then
$\gamma(s) \in B(x, r +\varepsilon)$ and
$$
|f(x)-f(x_0)| \leq (r +\varepsilon)
\sup_{y\in B(x, r +\varepsilon)} |(\nabla f)(y)|.$$
Taking $\varepsilon\rightarrow 0$ we get
$$
|f(x)-f(x_0)| \leq r\sup_{y\in B(x, r)}|(\nabla f)(y)|.$$
Since $f$ is smooth the supremum is finite which gives the first claim
of the lemma.  If $(\nabla f)(x_0) = 0$, then we can apply the
first part to $X_if$ and get
$$
\sup_{y\in B(x, r)}|(\nabla f)(y)| 
\leq Cr \sup_{y\in B(x, r)}|(\nabla\nabla f)(y)|,$$
$$
|f(x)-f(x_0)| \leq Cr^2\sup_{y\in B(x, r)}|(\nabla\nabla f)(y)|$$
which gives the second claim.
\qed

\begin{theorem}
\label{6.noLS}
Let $N$ be a stratified group and $\phi$ be a smooth homogeneous
norm on $N$.  For $\beta > 0,$ $p \geq 1$ put $\mu_{\beta,p} = \exp(-\beta \phi^p)/Z d\lambda$, where $Z$ is a normalizing factor such that $\mu_{\beta,p}$ is
a probability measure.  The measure $\mu_{\beta,p}$ satisfies no LS$_q$
inequality with $q\in(1,2]$.
\end{theorem}

{\it Proof}:  Fix $\beta>0$, $p \geq 1$, $q\in(1,2]$.  Suppose
that $\mu_{\beta,p}$ satisfies LS$_q$.  We are going to show that this
leads to contradiction.  Let
 $x_0$ be such that $(\nabla \phi)(x_0) = 0$.  For $t>0$
put $r = t^{(-p+1)/2}$ and $f = \max(\min((2 - d(x,tx_0))/r, 1), 0)$.
By homogeneity and Lemma \ref{hom-taylor}
we have $\phi(x) - \phi(tx_0) \leq C_1r^{2}$ on 
$B(tx_0, 2r) = \{x: d(x,tx_0) \leq 2r \}$, so 
$\phi(x)^p - \phi(tx_0)^p \leq C_2$.
Consequently the exponential factor in $\mu_{\beta,p}$ is comparable
to a constant on support of $f$.
Also $|\nabla f| \leq r^{-1}$ and
$$
\mu_{\beta,p}|f|^q \approx r^Q\exp(-\beta \phi(tx_0)^p),$$
$$
\log(\mu_{\beta,p}|f|^q) \approx -t^p,$$
$$
\mu_{\beta,p}|\nabla f| \approx r^{-q}r^Q\exp(-\beta \phi(tx_0)^p),$$
$$
\mu_{\beta,p}(|f|^q\log(|f|^q\mu_{\beta,p}|f|^q) \approx
\int_{B(tx_0, r)} |f|^qt^p d\mu_{\beta,p} \approx
t^{p}r^Q\exp(-\beta \phi(tx_0)^p).$$
Using LS$_q$ we get
$$
t^{p}r^Q\exp(-\beta \phi(tx_0)^p) \leq M r^{-q}r^Q\exp(-\beta \phi(tx_0)^p)$$
for large $t$, so
$$
t^{p} \leq M r^{-q} = Mt^{-q(-p+1)/2}$$
for large $t$, and $p \leq q(p-1)/2$.  Since $p\geq 1$ and $q \leq 2$,
this implies
$p \leq p - 1$ which is a contradiction.
\qed


\section{Log Sobolev Inequalities for Heat Kernel \\ on the Heisenberg Group. } 
\label{6.Heat Kernel on Heisenberg Group}

The heat kernels bound of the following form 
$$
\frac1{C|B(e,t^{1/2})|} e^{-\sigma {d^2(x)}{t}}\leq p(x,t) 
\leq \frac{C}{|B(e,t^{1/2})|} e^{-\frac1\sigma {d^2(x)}{t}}
$$
were well known since a few decades, see e.g. \cite{D}, \cite{VSCC} and references therein.
While the measures corresponding to the densities on the left and right have nice properties and in particular satisfy Poincar\'e and Logarithmic Sobolev inequality, this kind of sandwich bound does not imply
similar properties for the measure corresponding to the density in the middle.
Namely on a stratified groups one can write:
$$
C^{-1}p(x,t/\sigma) \leq \frac{1}{|B(e,t^{1/2})|}\exp(-\phi^2(x)/t)
\leq Cp(x,\sigma t)
$$
where $C$, $\sigma \geq 1$ are constants and $\phi$ is a smooth homogeneous
norm.  In Theorem \ref{6.noLS} we
proved that the density in the middle does not satisfy 
Logarithmic Sobolev inequality.  We give another example in the Appendix I.

In \cite{Li} it was observed that asymptotics from \cite{HuM} 
imply the following precise bound (extending \cite{BeGaGre}) on the heat
kernel $p$ (at time $t=1$) on the three-dimensional Heisenberg group $H_1$:

\begin{itemize}
 \item{(HK)}

{\textit{There exists a constant $L\in(0,\infty)$ such that for any} $x\equiv(\mathbf{x},z)\in H_1$ 
$$
L^{-1} \left(1+||\mathbf{x}||d(x) \right)^{-\frac12} e^{-\frac{d^2(x)}4}
\leq p(x) \leq L \left(1+||\mathbf{x}||d(x) \right)^{-\frac12} e^{-\frac{d^2(x)}4}
$$}

\end{itemize}

\noindent Let $d\nu_0\equiv \rho_0d\lambda \equiv e^{-\frac{d^2(x)}4} d\lambda /Z$ and set $d\mu = p d\lambda$.

\begin{theorem}\label{thm6.1}
There exist constants $C_1,C_2,D_1,D_2\in(0,\infty)$ such that
$$
\mu \left( f^2d^2\right)\leq C_2\mu |\nabla f|^2 + D_2\mu f^2
$$
and
$$
\mu \left( |f|d\right)\leq C_1\mu |\nabla f| + D_1\mu |f|
$$
\end{theorem}

\textit{Proof} :
Put $W = \frac{-1}{2}\log(1+\varepsilon||x||d)$
for some $ \varepsilon\in(0,1)$ to be chosen later.  We have
$$
|\nabla W|^2 = \varepsilon^2
\frac{|d\nabla||x|| + ||x||\nabla d|^2}{(1+\varepsilon||x||d)^2}$$
$$
\leq
\varepsilon^2
\frac{d^2 + ||x||^2}{(1+\varepsilon||x||d)^2} \leq 
\varepsilon^2d^2 + 1
$$
so, if $\varepsilon$ is small enough
$W$ satisfies assumptions of Theorem \ref{thm2.2}.

Now we observe that for $\varepsilon\in(0,1)$, we have
$$
(1+ ||x||d)^{-\frac12}\leq (1+\varepsilon||x||d)^{-\frac12}
\leq \frac1\varepsilon (1+ ||x||d)^{-\frac12}$$
This together with (HK) imply we can write $\mu = \exp(-W -V)\mu_0$
and apply Theorem \ref{thm2.2} to get the first claim.  We get the second
claim using Theorem \ref{thm2.6}.
\qed

By similar arguments as in Section 3 we obtain the following result

\begin{theorem}\label{thm6.P1}
Let $d\mu \equiv pd\lambda$.
There exist constants $M\in(0,\infty)$ such that
\begin{equation}\label{eqn6.5}
\phantom{AAAAAAAAA} M \mu(f-\mu f)^2 \leq \mu |\nabla f|^2 
\end{equation}
\end{theorem}

\smallskip
We are now  ready to prove the Log-Sobolev inequality for the heat kernel measure.
\begin{theorem}\label{thm6.2}
There exists a constant  $c\in(0,\infty)$ such that on Heisenberg group $H_n$
we have
$$
\mu \left( f^2\log \frac{f^2}{\mu f^2}\right)\leq c\mu |\nabla f|^2 
$$
\end{theorem}

\noindent\textbf{Remark}: The case of $H_1$ is proven in \cite{Li}.  While
our proof uses heat kernel estimates from \cite{Li}, in \cite{Li} large
part is devoted to proof of estimate (\ref{2.z}) for heat kernel
measure on $H_1$ -- using our methods we could give different proof
for this part, but instead
we work directly with Log-Sobolev inequality.

\noindent\textit{Proof}:  First consider $H_1$.  In the proof of Theorem \ref{thm6.1} we
wrote 
$\mu = e^{-W -V}\mu_0$.  Consider now $\mu_1 = e^{-W}\mu = e^{-U}d\lambda$.
$\mu_1$ satisfies Log-Sobolev inequality as a consequence of Theorem \ref{thm4.1}.
The result for $H_1$ follows, since $\mu$ is equivalent to $\mu_1$.

Now, write $H_n = G/N$, where $G=\prod_{i=1}^nH_1$,
$N=\{((0,z_1),\dots,(0,z_n)): \sum z_i = 0 \}$ and let $\pi$ be the
canonical homomorfizm from $G$ to $H_n$.  Since heat
kernel on $H_n$ is an image of product of heat kernels
on $G = \prod_{i=1}^nH_1$, and since Log-Sobolev inequality holds on
product, we have
$$
\mu_{H_n} \left( f^2\log \frac{f^2}{\mu_{H_n} f^2}\right) =
\mu_G \left( (f\circ \pi)^2\log \frac{(f\circ \pi)^2}{\mu_G(f\circ \pi)^2}
\right)
$$
$$
\leq c\mu_G |\nabla (f\circ \pi)|^2 = c\mu_G |(\nabla f)\circ \pi|^2 =
c\mu_{H_n}|\nabla f|^2.$$

\qed

\section{Appendix:  Examples of No Spectral Gap. } 
\label{7.No Spectral Gap}
In case of measures on real line the following necessary and sufficient condition
for Poincar\'e inequality 
to hold was provided by Muckenhoupt \cite{Mu} (\cite{ABC00})
which in the special case of a measure $d\mu\equiv \rho dx $ 
 can be stated as follows: Given $q\in[1,\infty)$ and $\frac1q+\frac1p=1$
\begin{equation}\label{A1.1}
\exists C\in(0,\infty)\ \ \ \mu|f-\mu f|^q \leq \mu |f'|^q 
\Longleftrightarrow  
B_\pm\equiv\sup_{r\in\mathbb{R}^\pm} B_\pm(r)
\end{equation}
where
$$
B_\pm(r) \equiv
\left(\mu([r,\pm\infty))\right)^\frac1q \cdot \left(\int_{[0,\pm r]}\rho^{-\frac{p}{q}}\right)^\frac1p < \infty
$$
Consider $\rho\equiv e^{-U}dx/Z$ with $U\equiv \beta|x|^p(1 +\varepsilon\cos x)$, defined $\varepsilon\in(0,1)$ and some $\beta\in(0,\infty)$.
Then, with $r = 2n\pi +\frac\pi2$, we have
$$
B_+(r) \ > \ \left(\int_{2n\pi +\frac43\pi}^{2n\pi + \frac83\pi} 
e^{-\beta |x|^p(1-\frac\varepsilon2)} dx \right)^\frac1q  \cdot
\left(\int_{2n\pi -\frac23\pi}^{2n\pi + \frac23\pi}
e^{+\frac{p}{q}\beta |x|^p(1+\frac\varepsilon2)} dx \right)^\frac1p 
$$
$$
>
e^{-\beta\frac1q |2n\pi + \frac83\pi|^p(1-\frac\varepsilon2)} \left(\frac43\pi\right)^\frac1q
\cdot
e^{+\frac{1}{q}\beta |2n\pi - \frac23\pi|^p(1+\frac\varepsilon2)}
\left(\frac43\pi\right)^\frac1p
$$
$$
=\frac43\pi \exp\left\{ \frac{\beta(2n\pi)^p}{q} 
\left[ |1 - \frac1{3n}|^p(1+\frac\varepsilon2)
-|1 + \frac4{3n}|^p(1-\frac\varepsilon2)
\right]\right\}
$$
$$
\sim \frac43\pi \exp\left\{ \frac{\beta(2n\pi)^p}{q} (\varepsilon + o(\frac1n))\right\} \to \infty \qquad  as \qquad n\to\infty
$$
Alternatively one can study lower bound asymptotic for $B_\pm$ thinking of
$U=V+\delta V$ as a perturbation of $V\equiv\beta|x|^p $ as follows.
We notice that by Jensen inequality
$$
B_+(r,U) \geq B_+(r,V)\exp\left\{ - \frac1q\beta\frac{\int_r^\infty \delta V e^{-V}dx}{\int_r^\infty e^{-V}dx} 
+ \frac{p}{q}\beta\varepsilon\frac{\int_0^r \delta Ve^{+V}dx}{\int_0^r e^{+V}dx}  \right\}
$$
Hence one can use a procedure based essentially on integration by parts 
to study the integrals in the exponential.
For example in case $p=2$ one gets the following an asymptotic lower bound
$$
B_+(r,U)\geq B_+(r,V)\exp\{-\beta\varepsilon r\cos r + O(1)\}
$$
We summarise our considerations in the above as follows
\begin{proposition}\label{proA1}
Suppose $p\geq 1$. In any neighbourhood
$$
\frac1C e^{-{(1+\delta)}\beta|x|^p}\leq \rho \leq 
C e^{-\frac1{1+\delta}\beta|x|^p}
$$
with  arbitrary $\delta\in(0,1)$ and  some $C\in(1,\infty)$, of a measure $ d\mu_0\equiv \frac{e^{-\beta|x|^p}dx}{Z}$ satisfying the Poincar\'e inequality 
there is a measure $d\mu\equiv \rho dx$ for which this inequality fails. 
\hfill{$\circ$}
\end{proposition}
The example provided above illustrates similar phenomenon for other coercive inequalities.


 
\end{document}